\documentclass[12pt]{article}

\usepackage{amsmath}
\usepackage{amssymb}
\usepackage{setspace}
\usepackage{ifthen}
\usepackage{bm}
\usepackage{bbm}
\usepackage{url}
\usepackage{multirow}

\usepackage{epsfig}
\usepackage{epstopdf}
\usepackage{graphicx}
\usepackage{multirow}
\usepackage{verbatim}
\usepackage{natbib}
\usepackage{enumerate}

\usepackage{url}
\usepackage{calligra}
\usepackage{algorithm}
\usepackage{algpseudocode}
\usepackage{multicol}
\usepackage{subfig}

\usepackage[table]{xcolor}
\usepackage{tikz}
\usepackage{rotating}

\newcommand{\be}{\begin{eqnarray}}
\newcommand{\ee}{\end{eqnarray}}
\newcommand{\ba}{\begin{eqnarray*}}
\newcommand{\ea}{\end{eqnarray*}}

\newcommand{\gcl}{\cellcolor{gray!75}}
\newcommand{\dcl}{\cellcolor{gray!25}}

\algrenewcommand\alglinenumber[1]{\scriptsize #1:} 

\begin{document}


\title{
\begin{center}
{
\textbf{\LARGE{Demand Estimation from Sales Transaction Data - Practical Extensions}}\\
\vspace{.5cm}
\vspace{.5cm}
\textit{Submission to Journal of Revenue and Pricing Management}
\vspace{.5cm}
}
\end{center}
}

\date{}

\author{
Norbert Remenyi\\
Xiaodong Luo\\}

\maketitle{
\begin{center}
Sabre\\
\vspace{.5cm}
May 31, 2020\\
\end{center}
}

\begin{abstract}
\noindent In this paper we discuss practical limitations of the standard choice-based demand models used in the literature to estimate demand from sales transaction data. We present modifications and extensions of the models and discuss data preprocessing and solution techniques which are useful for practitioners dealing with sales transaction data. Among these, we present an algorithm to split sales transaction data observed under partial availability, we extend a popular Expectation Maximization (EM) algorithm for non-homogeneous product sets, and we develop two iterative optimization algorithms which can handle much of the extensions discussed in the paper.
\end{abstract}
\noindent{\em Keywords.} Demand estimation. Demand untruncation. Multinomial logit model. EM algorithm. MM algorithm. Frank-Wolfe method. Revenue management

\section{Introduction}
\noindent Demand estimation using censored sales transaction data has many applications in airline commercial planning process. We refer readers to a survey paper by \citet{Azadeh2014} for a brief introduction to this topic. In this paper we discuss some practical limitations and extensions of a particular choice-based demand model popular in the literature to estimate demand from sales transaction data. Discrete choice models (e.g., \citet{benakiva1994}, \citet{train2003}) have provided a popular approach for estimating demand for different products in a set of substitutable items, especially in transportation and revenue management applications. We look at a common demand model, which appeared in several papers, including \citet{dai2014}, \citet{vvrr12}, and \citet{tarek2016}. The motivation of this work came from observing a few shortcoming of these models when applied in practice, on airline revenue management data.

We will build on and extend the work presented in \citet{vvrr12} and \citet{tarek2016}. They combine a multinomial logit (MNL) choice model with non-homogeneous Poisson arrivals over multiple periods. The MNL model has been used by many practitioners and researchers to represent the underlying choice behavior. Although its property of independence of irrelevant alternatives (IIA) is somewhat restrictive, the model is simple, leading to tractable estimation and assortment optimization \citep{vanryzin04}. The problem is how to jointly estimate the arrival rates of customers and the preference weights of the products via maximizing the likelihood functions. There are two different likelihood functions. The first one is the incomplete data likelihood function (see $(\ref{eq:incomp})$ below and also Eq. (2) in \cite{vvrr12}) and the second one is the log-likelihood function which is based on the primary  demand (see $(\ref{eq:compLL})$ below and also Eq. (13) in \cite{vvrr12}). The inputs are observed historical sales, availability of the products, and market share information.

Our contribution is to discuss practical limitations of the specific model above, and present some interesting extensions. We will discuss partial availability of products, some relaxation of the IIA assumption, constrained parameter space, non-homogeneous product set, and the interpretation of the no-purchase option and related market share. We hope this discussion can facilitate more research and extension of these models.

We will present an algorithm to split sales transaction data observed under partial availability, and an extension of the EM algorithm for the case when we observe a non-homogeneous product set. We develop two iterative optimization algorithms which incorporate partial availability information, non-homogeneous product set, ability to control the availability of outside alternative, and an upper bound on the arrival rates of customers. In the first formulation we use a market share constraint at each time period, and incorporate them into the objective function through the preference weights of the outside alternative. The formulation is solved using the Frank-Wolfe algorithm, leading to a simple coordinate descent algorithm. In the second formulation we use a single, aggregate market share constraint over the time horizon, and assume knowledge of preference weights of the outside alternative. Using this formulation we develop a fast, iterative minorization-maximization algorithm (MM) building on the work in \citet{tarek2016}.

While the EM algorithm focuses on solving the complete data log-likelihood functions,  the two new algorithms (like \citet{tarek2016}) aim to solve the incomplete likelihood function directly. We remark that both likelihood functions render similar quality solutions, but they involve different intermediate decision variables and require different solution approaches. It is not known that these two methods are mathematically equivalent.

\section{Practical limitations of existing models}
In this section we discuss in detail some practical limitations of the choice-based demand model discussed in \citet{vvrr12} and in \citet{tarek2016}. The model combines non-homogeneous Poisson arrivals over multiple periods with a multinomial logit (MNL) choice model. The model assumes a retailer offers a fixed number of $n$ products over a time horizon $T$, and the products are either available or not available for sale in a time period. From the observed sales data we estimate the demand, the sales we would have observed if all the products were available for purchase. The total demand of all the products at time $t$ (including outside alternatives and the no-purchase option) is modeled as a Poisson distributed random variable with parameter $\lambda_t$. Hence we model the total demand as a non-homogeneous Poisson model, since different time periods are allowed to have different mean demands.

The demand for individual products is modeled using a multinomial distribution, given the total demand of all products. The preferences for different products are assumed to be fixed over time, and the probability that customer chooses product $i$ when $i \in S_t$ is modeled through the simple multinomial logit (MNL) model, that is
\be
\label{eq:BAM}
P_i(S_t, \bm v) = \frac{v_i}{v_0 + \sum_{j \in S_t} v_j}
\ee
When $i \notin S_t$, then $P_i(S_t, \bm v) = 0$. Here $v_i,~i=1,...,n$ is the positive preference weight for product $i$, and $S_t$ is the set of products available for sale at time $t$. The preference weight of outside alternatives and no-purchase option (OA) are embedded into the coefficient $v_0$. The parameter is set to $v_0=1$ in the references above, following a standard approach of normalizing.

The incomplete data likelihood function of the model is defined as
\footnotesize
\ba
L_I(\bm v, \bm \lambda) &=& \prod_{t=1}^{T} \left[ \left( P(m_t~\textnormal{customers buy in period}~t|\bm v, \bm \lambda)
\frac{m_t!}{z_{1t}!z_{2t}!\cdots z_{nt}!} \right) \prod_{j \in S_t} \left[ \frac{P_{j}(S_t, \bm v)}{\sum_{i \in S_t} P_{i}(S_t, \bm v)} \right]^{z_{jt}} \right]
\ea
\normalsize
where
\footnotesize
\ba
P(m_t~\textnormal{customers buy in period}~t|\bm v, \bm \lambda) &=& \frac{\left[ \lambda_t \sum_{i \in S_t} P_i(S_t, \bm v) \right]^{m_t}
\exp\left(- \lambda_t \sum_{i \in S_t} P_i(S_t, \bm v) \right)}{m_t!}
\ea
\normalsize
In the above equations $z_{it}$ denotes the number of purchases of product $i$ at time period $t$, and $m_t = \sum_{i=1}^n z_{it}$ denotes the total number of purchases in period $t$. After some algebra the log-likelihood function can be written as
\footnotesize
\be
\label{eq:loglik}
l_I(\bm{v}, \bm{\lambda}) = \sum_{t=1}^{T} \left[ m_t \log\left( \frac{\lambda_t}{v_{0} + \sum_{i \in S_t} v_i} \right) - \lambda_t \frac{\sum_{i \in S_t} v_i}{v_{0} + \sum_{i \in S_t} v_i} +\sum_{i \in S_t} z_{it}\log(v_i) \right]
\ee
\normalsize

One approach is to directly maximize the log-likelihood function and jointly estimate the preference weights of the products and the arrival rates of customers. The above log-likelihood function, however, is hard to solve in general, therefore the research literature discusses different approaches to estimate the parameters of this model, given sales data and information on what was available for sale. \citet{vvrr12} developed an elegant EM algorithm by looking at the problem in terms of primary (first choice) demand, and treating the observed demand as incomplete observations of primary demand. \citet{tarek2016} solves the estimation problem by specializing the minorization-maximization (MM) procedure, which is an iterative algorithm for maximizing an objective function by successively maximizing a simpler function that minorizes the true objective function.

It is also interesting  to note that the objective function has a continuum of maximizers. For this reason, \citet{vvrr12} and \citet{tarek2016} imposed additional constraint on the preference weights as a function of the market share
\ba
s = \frac{\sum_{i=1}^n v_i}{v_0 + \sum_{i=1}^n v_i}
\ea

There are a number of limiting assumptions and possible extensions of the model presented above, when we apply it to real data.

\begin{enumerate}
\item \textbf{Products are fully open or closed for sale}
\begin{itemize}
\item[] The discussed model assumes that a product is either available or not for sale. Practitioners often work with aggregated data, and unable to capture the sales at a very granular level, every time the assortment changes. It is very practical to extend the model to consume data where we observe partial availability. For example, a product can be 80\% open for sale in a time period.
\end{itemize}
\item \textbf{MNL assumption}
\begin{itemize}
\item[] Customers choose from the available products according to an MNL model. One of the properties of the MNL model is the independence of irrelevant alternatives (IIA), which can be unrealistic in real applications. If customers, for instance, always purchase the product with the lowest price, the algorithm in \citet{vvrr12} does not converge. If the sell-down is strong, the demand is grossly overestimated.
\end{itemize}
\item \textbf{Unconstrained parameter space}
\begin{itemize}
\item[] In practice we can encounter parameter estimates which are unreasonable in a business scenario. This may be due to the fact that the underlying assumptions of the model (such as MNL) do not exactly fit the true data generating model. Instead of using a more sophisticated model and develop its solution algorithm, we might want to simply constrain the parameters of the simpler model. \citet{tarek2016} discusses regularization as an elegant solution to the problem, and they develop algorithms for $L_1$ regularization (Lasso regression) and $L_2$ regularization (Ridge regression). It can be also of interest to solve the problem by putting an upper bound on the arrival rate of customers ($\lambda_t$). This can be helpful to regularize the model by having an interpretable business fence on these parameters.
\end{itemize}
\item \textbf{Homogeneous product set}
\begin{itemize}
\item[] The model assumes that over different time periods we have the same set of existing products. Some might be unavailable for sale, but there exist an underlying demand for them. In revenue management, we often encounter changing product IDs (unique flight identifiers) due to schedule changes, hence we need to be able to model a non-homogeneous product set over time. Instead of handling the homogeneous parts separately, we would like to use the data over a larger time interval to borrow power to estimate the parameters.
\end{itemize}
\item \textbf{No-purchase option is always available}
\begin{itemize}
\item[] The model above assumes that the no-purchase option is always available, that is fully open. If we include competitor's products into the no-purchase option, this can be an unrealistic assumption. For instance, airline competitors likely control their inventory similarly as the host airline, gradually closing down products as we get closer to departure. This assumption has implications on how the host market share is interpreted in the models.
\end{itemize}
\end{enumerate}

In this paper we address some of these points, and discuss how we could relax these assumptions, incorporate them into the model, and estimate the parameters. These ideas might foster further extensions of the existing models and rigorous research in the future.

\subsection{Partial availability}
Practitioners can encounter data at an aggregate level, where we can observe partial availability of products. As an example, in a revenue management system we store information on sales and availability at certain pre-departure time points. It is possible that a product becomes closed for sale during the time period, and we observe partial availability. For instance, a booking class on a flight can be open 60\% of the time in a time interval, and the bookings we observe are matched to this partial availability. It would be of interest to extend the algorithms to work with this type of data and handle the full spectrum of availability in $[0,1]$, as opposed to be restricted to either open (1) or closed (0) products. Another possible venue is to modify the data to fit the algorithms already developed in the literature.

\subsubsection{Extending the attraction model}
\label{sec:extend_attr_model}

A natural way to incorporate partial availability is to extend the MNL purchase probabilities as
\ba
P^{\star}_{jt}(S_t, \bm v, \bm o) = \frac{v_j \cdot o_{jt}}{v_0 + \sum_{i \in S_t} v_i \cdot o_{it}}
\ea
where $o_{it} \in \left[0,1\right]$ represents availability of product $i$ at time $t$. It is obvious that $o_{it} > 0$ is equivalent to $i \in S_t$. This simple formulation modifies the purchase probabilities by linearly adjusting preference weight $v_i$ with open percentage $o_{it}$. If the open percentage is zero or one, we get back to the original formulation.

We have not derived the EM algorithm of \citet{vvrr12} with the extended purchase probability definition, but this could be an interesting venue for research. For demonstration purposes we can resort to directly maximizing the log-likelihood using a solver. The modified likelihood function becomes

\scriptsize
\be
\label{eq:incomp}
L_I(\bm v, \bm \lambda) &=& \prod_{t=1}^{T} \left[ \left( P(m_t~\textnormal{customers buy in period}~t|\bm v, \bm \lambda)
\frac{m_t!}{z_{1t}!z_{2t}!\cdots z_{nt}!} \right) \prod_{j \in S_t} \left[ \frac{P^{\star}_{j}(S_t, \bm v, \bm o)}{\sum_{i \in S_t} P^{\star}_{i}(S_t, \bm v, \bm o)} \right]^{z_{jt}} \right]
\ee
\normalsize
where
\footnotesize
\ba
P(m_t~\textnormal{customers buy in period}~t|\bm v, \bm \lambda) &=& \frac{\left[ \lambda_t \sum_{i \in S_t} P^{\star}_{i}(S_t, \bm v, \bm o) \right]^{m_t}
\exp\left(- \lambda_t \sum_{i \in S_t} P^{\star}_{i}(S_t, \bm v, \bm o) \right)}{m_t!}
\ea
\normalsize

\noindent and, after omitting the constant terms, the log-likelihood function modifies to
\footnotesize
\be
\label{eq:loglik_openpct}
l_I(\bm{v}, \bm{\lambda}) = \sum_{t=1}^{T} \left[ m_t \log\left( \frac{\lambda_t}{v_{0} + \sum_{i \in S_t} v_i \cdot o_{it}} \right) - \lambda_t \frac{\sum_{i \in S_t} v_i \cdot o_{it}}{v_{0} + \sum_{i \in S_t} v_i \cdot o_{it}} + \sum_{i \in S_t} z_{it}\log(v_i) \right]
\ee
\normalsize

\subsubsection{Data disaggregation}
\label{sec:data_disagg}

In the previous section we discussed a simple formulation which could potentially be used to incorporate partial availability information into the purchase probability definition. Another approach to handle partial availability is to disaggregate the data to fully open and closed assortments, and use existing algorithms for the estimation. We can split the observed sales under partial availability by making a simple assumption that the sales are distributed uniformly over time.

To demonstrate this idea on a simple example, let us assume that the observed sales $b$ and open percentages $o$ for three available products are
\small
\ba
\bm b = \begin{bmatrix}
 1 \\ 2 \\ 5
\end{bmatrix}, \quad
\bm o = \begin{bmatrix}
 1.0 \\ 0.8 \\ 0.5
\end{bmatrix}
\ea
\normalsize
We assumed that the elements of $\bm o$ are non-increasing from top to bottom. Then we can represent $\bm o$ as the sum of three fully open and closed assortments with weights $\alpha_i$ as
\small
\ba
\bm o = \begin{bmatrix}
 1.0 \\ 0.8 \\ 0.5
\end{bmatrix}
= \alpha_1  \begin{bmatrix}
 1 \\ 0 \\ 0
\end{bmatrix}
+ \alpha_2  \begin{bmatrix}
 1 \\ 1 \\ 0
\end{bmatrix}
+ \alpha_3  \begin{bmatrix}
 1 \\ 1 \\ 1
\end{bmatrix}
= 0.2  \begin{bmatrix}
 1 \\ 0 \\ 0
\end{bmatrix}
+ 0.3  \begin{bmatrix}
 1 \\ 1 \\ 0
\end{bmatrix}
+ 0.5  \begin{bmatrix}
 1 \\ 1 \\ 1
\end{bmatrix}
\ea
\normalsize

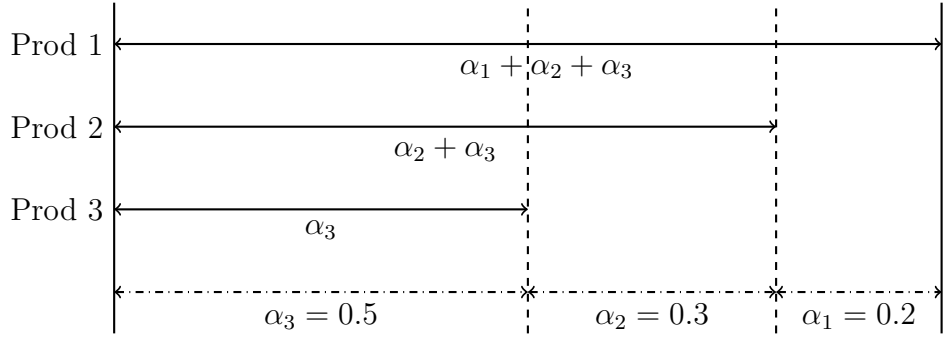
\begin{figure}
\centering
\begin{tikzpicture}[scale = 0.55]
\draw[step=0.5cm,white,very thin] (0,1) grid (20,9);
\draw[thick,-] (0,1) -- (0,9);
\draw[thick,-] (20,1) -- (20,9);
\draw[thick,<->] (0,8) -- node[below, xshift = 0.25cm] {$\alpha_1 + \alpha_2 + \alpha_3$} (20,8);
\draw[thick,<->] (0,6) -- node[below] {$\alpha_2 + \alpha_3$} (16,6);
\draw[thick,<->] (0,4) -- node[below] {$\alpha_3$} (10,4);
\draw[dashdotted,black,thick,<->] (0,2) -- node[below, text = black] {$\alpha_3 = 0.5$} (10,2);
\draw[dashdotted,thick,<->] (10,2) -- node[below, text = black] {$\alpha_2 = 0.3$} (16,2);
\draw[dashdotted,thick,<->] (16,2) -- node[below, text = black] {$\alpha_1 = 0.2$} (20,2);
\draw[black,thick,dashed] (10,1) -- (10,9);
\draw[black,thick,dashed] (16,1) -- (16,9);
\draw (0.01pt,4) -- (-0.01pt,4) node[anchor=east] {Prod $3$};
\draw (0.01pt,6) -- (-0.01pt,6) node[anchor=east] {Prod $2$};
\draw (0.01pt,8) -- (-0.01pt,8) node[anchor=east] {Prod $1$};
\end{tikzpicture}
\caption{Split example}
\label{fig:Split_ex}
\end{figure}

Note that $\alpha_i,~i=1,...,n$ represent time proportions and $\sum_{i=1}^n \alpha_i = 1$. The results indicate that all products were available for sale 50\% of the time, two products were available 30\% of the time, and one product was available 20\% of the time. A graphical representation of this example is shown in Figure \ref{fig:Split_ex}.

For the general case, the elements of $\bm \alpha$ can be calculated as the consecutive differences in the open percentages, that is
\ba
\alpha_i &=& o_i - o_{i+1}, \quad i = 1,\ldots,n-1 \\
\alpha_n &=& o_n
\ea
\noindent Following this simple idea we can split the observed sales under partial availability using the following identity
\small
\ba
\bm b &=& \begin{bmatrix}
 \frac{\alpha_1}{\sum_{i=1}^n \alpha_i} \\ 0 \\ \vdots \\ 0 \\ 0
\end{bmatrix}
\circ \begin{bmatrix}
 b_1 \\ b_2 \\ \vdots \\ b_{n-1} \\ b_n
\end{bmatrix}
+ \begin{bmatrix}
 \frac{\alpha_2}{\sum_{i=1}^n \alpha_i} \\ \frac{\alpha_2}{\sum_{i=2}^n \alpha_i} \\ \vdots \\ 0 \\ 0
\end{bmatrix}
\circ \begin{bmatrix}
 b_1 \\ b_2 \\ \vdots \\ b_{n-1} \\ b_n
\end{bmatrix}
+ \ldots
+ \begin{bmatrix}
 \frac{\alpha_{n-1}}{\sum_{i=1}^n \alpha_i} \\ \frac{\alpha_{n-1}}{\sum_{i=2}^n \alpha_i} \\ \vdots \\ \frac{\alpha_{n-1}}{\sum_{i=n-1}^n \alpha_i} \\ 0
\end{bmatrix}
\circ \begin{bmatrix}
 b_1 \\ b_2 \\ \vdots \\ b_{n-1} \\ b_n
\end{bmatrix}
+ \\
&& \begin{bmatrix}
 \frac{\alpha_n}{\sum_{i=1}^n \alpha_i} \\ \frac{\alpha_n}{\sum_{i=2}^n \alpha_i} \\ \vdots \\ \frac{\alpha_n}{\sum_{i=n-1}^n \alpha_i} \\ \frac{\alpha_n}{\alpha_n}
\end{bmatrix}
\circ \begin{bmatrix}
 b_1 \\ b_2 \\ \vdots \\ b_{n-1} \\ b_n
\end{bmatrix}
\ea
\normalsize
where $\circ$ denotes the element-wise multiplication, or Hadamard product. Note that if $\alpha_i=0$, we do not need a split to create a new assortment. The sales splitting algorithm under partial availability is described in Algorithm \ref{algo_sales_split}.

\begin{algorithm}[!htbp]
\footnotesize 
\caption{Sales splitting algorithm for partial availability}
\label{algo_sales_split}
\begin{algorithmic}[1]
\State $\bm b$: vector of observed sales
\State $\bm o$: vector of observed open percentages
\State Sort $\bm o$ in decreasing order and apply the order to $\bm b$
\State \textbf{Calculate time proportion vector $\bm \alpha$}
\For{$i = 1,\ldots,n$}
    \If{$i = n$}
        \State $\alpha_i = o_i$
    \Else
        \State $\alpha_i = o_i - o_{i+1}$
    \EndIf
\EndFor
\State \textbf{Split sales}
\For{$i = 1,\ldots,n$}
    \If{$\alpha_i = 0$}
        \State Continue
    \Else
        \State Calculate partial sales vector $\bm{b^{(i)}}$ and open percentage vector $\bm{o^{(i)}}$
         \For{$j = 1,\ldots,n$}
            \If{$j \leq i$}
            \State $b^{(i)}_j = \frac{\alpha_i}{o_j} b_j$
            \State $o^{(i)}_j = 1$
            \Else
            \State $b^{(i)}_j = 0$
            \State $o^{(i)}_j = 0$
            \EndIf
         \EndFor
    \EndIf
\EndFor
\end{algorithmic}
\end{algorithm}

The algorithm results in partial sales vectors $\bm{b^{(i)}}$ and fully open and closed assortment vectors $\bm{o^{(i)}}$ for $i \in \{j:\alpha_j > 0\}$. The splitting logic ensures that
\ba
\bm{b} = \sum_{i \in \{j:\alpha_j > 0\}} \bm{b^{(i)}}
\ea
 After the splitting logic, we could combine intervals from different time $t$ which have exactly the same product offerings and product availabilities, to improve performance.

\subsubsection{Comparison on a simulated example}

To demonstrate the estimation from sales data with partial availability, we extended the example from \citet{vvrr12} by adding partial availability data to the observed sales. The observed data is presented in Table \ref{ex:vvrr_op}.

\begin{table}[!htbp]
\centering
\resizebox{0.8\linewidth}{!}{
\begin{tabular}[h]{|c||c|c|c|c|c|c|c|c|c|c|c|c|c|c|c|}
	\hline
\multirow{2}{*}{\bf Sales} & \multicolumn{15}{c|}{\bf Period} \\ \cline{2-16}
            & \bf 15 & \bf 14 & \bf 13 & \bf 12 & \bf 11 & \bf 10 & \bf 9 & \bf 8 & \bf 7 & \bf 6 & \bf 5 & \bf 4 & \bf 3 & \bf 2 & \bf 1 \\
\hline
1  & 10 & 15 & 11 & 14 & 0  & 0  & 0  & 0  & 0  & 0  & 0  & 0  & 0  & 0  & 0 \\
2  & 11 & 6  & 11 & 8  & 20 & 16 & 0  & 0  & 0  & 0  & 0  & 0  & 0  & 0  & 0 \\
3  & 5  & 6  & 1  & 11 & 4  & 5  & 14 & 7  & 11 & 0  & 0  & 0  & 0  & 0  & 0 \\
4  & 4  & 4  & 4  & 1  & 6  & 4  & 3  & 5  & 9  & 9  & 6  & 9  & 0  & 0  & 0 \\
5  & 0  & 2  & 0  & 0  & 1  & 0  & 1  & 3  & 0  & 3  & 3  & 5  & 2  & 3  & 3  \\ \hline
\bf Open percentage & \bf 15 & \bf 14 & \bf 13 & \bf 12 & \bf 11 & \bf 10 & \bf 9 & \bf 8 & \bf 7 & \bf 6 & \bf 5 & \bf 4 & \bf 3 & \bf 2 & \bf 1 \\
\hline
1  & 0.7 & 0.3 & 1.0 & 1.0 & 0.0 & 0.0 & 0.0 & 0.0 & 0.0 & 0.0 & 0.0 & 0.0 & 0.0 & 0.0 & 0.0 \\
2  & 0.8 & 0.5 & 1.0 & 1.0 & 0.3 & 0.3 & 0.0 & 0.0 & 0.0 & 0.0 & 0.0 & 0.0 & 0.0 & 0.0 & 0.0 \\
3  & 0.9 & 0.9 & 1.0 & 1.0 & 0.3 & 0.5 & 0.9 & 0.7 & 0.7 & 0.0 & 0.0 & 0.0 & 0.2 & 0.2 & 0.0 \\
4  & 1.0 & 0.9 & 1.0 & 1.0 & 0.5 & 1.0 & 1.0 & 0.8 & 0.9 & 0.6 & 0.5 & 0.2 & 0.3 & 0.5 & 0.0 \\
5  & 1.0 & 1.0 & 1.0 & 1.0 & 1.0 & 1.0 & 1.0 & 0.9 & 1.0 & 1.0 & 1.0 & 1.0 & 1.0 & 1.0 & 1.0  \\ \hline
\end{tabular}
}
\caption{Example with partial availability}
\label{ex:vvrr_op}
\end{table}

We first use the extended attraction model discussed in Section \ref{sec:extend_attr_model}. Direct maximization of the log-likelihood (\ref{eq:loglik_openpct}) by a nonlinear solver results in total primary demand of 1194.6. The detailed estimated demands and parameters are presented in Table \ref{vvrr_op_sol}.

\begin{table}[!htbp]
\centering
\resizebox{0.95\linewidth}{!}{
\begin{tabular}[!htbp]{|c||c|c|c|c|c|c|c|c|c|c|c|c|c|c|c||c|}
	\hline
\multirow{2}{*}{\bf Estimates} & \multicolumn{15}{c||}{\bf Period} & \multirow{2}{*}{$\mathbf{v}_i$} \\ \cline{2-16}
            & \bf 15 & \bf 14 & \bf 13 & \bf 12 & \bf 11 & \bf 10 & \bf 9 & \bf 8 & \bf 7 & \bf 6 & \bf 5 & \bf 4 & \bf 3 & \bf 2 & \bf 1 & \\ \hline
1 & 15.02 & 20.07 & 12.47 & 15.70 & 33.60 & 22.73 & 19.60 & 19.33 & 24.84 & 38.14 & 32.27 & 84.35 & 5.74 & 7.22 & 35.06 & 1.000 \\
2 & 11.23 & 15.01 & 9.32 & 11.74 & 25.13 & 17.00 & 14.66 & 14.45 & 18.58 & 28.53 & 24.14 & 63.08 & 4.29 & 5.40 & 26.22 & 0.748 \\
3 & 3.91 & 5.22 & 3.24 & 4.09 & 8.74 & 5.92 & 5.10 & 5.03 & 6.46 & 9.93 & 8.40 & 21.95 & 1.49 & 1.88 & 9.13 & 0.260 \\
4 & 1.97 & 2.63 & 1.63 & 2.06 & 4.41 & 2.98 & 2.57 & 2.53 & 3.26 & 5.00 & 4.23 & 11.06 & 0.75 & 0.95 & 4.60 & 0.131 \\
5 & 0.40 & 0.53 & 0.33 & 0.41 & 0.89 & 0.60 & 0.52 & 0.51 & 0.66 & 1.01 & 0.85 & 2.23 & 0.15 & 0.19 & 0.93 & 0.026 \\ \hline
$\mathbf{\lambda}_t$ & 46.48 & 62.10 & 38.57 & 48.57 & 103.95 & 70.32 & 60.65 & 59.79 & 76.84 & 118.01 & 99.84 & 260.94 & 17.76 & 22.34 & 108.48 & \\ \hline
\end{tabular}
}
\caption{Estimated demand and parameters using the extended attraction model}
\label{vvrr_op_sol}
\end{table}

Second, we demonstrate the data disaggregation procedure from Section \ref{sec:data_disagg}. We first split the observed sales data with partial availability to fully open and closed assortments, apply the EM algorithm of \citet{vvrr12}, and then aggregate the solution back. To demonstrate this, let us apply Algorithm \ref{algo_sales_split} on the data presented in Table \ref{ex:vvrr_op}. The disaggregated sales are shown in Table \ref{ex:vvrr_op_disagg}.

\begin{table}[!htbp]
\centering
\resizebox{0.90\linewidth}{!}{
\begin{tabular}[!htbp]{|c||cccc|cccc|c|c|ccc|ccc|cc|}
	\hline
\multirow{2}{*}{\bf Sales} & \multicolumn{18}{c|}{\bf Period} \\ \cline{2-19}
& \multicolumn{4}{c|}{\bf 15} & \multicolumn{4}{c|}{\bf 14} & \multicolumn{1}{c|}{\bf 13} & \multicolumn{1}{c|}{\bf 12} & \multicolumn{3}{c|}{\bf 11} & \multicolumn{3}{c|}{\bf 10} & \multicolumn{2}{c|}{\bf 9} \\
\hline
1 & \gcl & \gcl & \gcl & 10.00 & \gcl & \gcl & \gcl & 15.00 & 11.00 & 14.00 & \gcl & \gcl & \gcl  & \gcl & \gcl & \gcl  & \gcl & \gcl  \\
2 & \gcl & \gcl & 1.38 & 9.62  & \gcl & \gcl & 2.40 & 3.60  & 11.00 & 8.00  & \gcl & \gcl & 20.00 & \gcl & \gcl & 16.00 & \gcl & \gcl  \\
3 & \gcl & 0.56 & 0.56 & 3.89  & \gcl & 2.67 & 1.33 & 2.00  & 1.00  & 11.00 & \gcl & \gcl & 4.00  & \gcl & 2.00 & 3.00  & \gcl & 14.00 \\
4 & 0.40 & 0.40 & 0.40 & 2.80  & \gcl & 1.78 & 0.89 & 1.33  & 4.00  & 1.00  & \gcl & 2.40 & 3.60  & 2.00 & 0.80 & 1.20  & 0.30 & 2.70  \\
5 & 0.00 & 0.00 & 0.00 & 0.00  & 0.20 & 0.80 & 0.40 & 0.60  & 0.00  & 0.00  & 0.50 & 0.20 & 0.30  & 0.00 & 0.00 & 0.00  & 0.10 & 0.90  \\ \hline
\end{tabular}}
\resizebox{0.90\linewidth}{!}{
\begin{tabular}[!htbp]{|c||ccc|ccc|cc|cc|cc|ccc|ccc|c|}
\hline
\bf \phantom{Sales} & \multicolumn{3}{c|}{\bf 8} & \multicolumn{3}{c|}{\bf 7} & \multicolumn{2}{c|}{\bf 6} & \multicolumn{2}{c|}{\bf 5} & \multicolumn{2}{c|}{\bf 4} & \multicolumn{3}{c|}{\bf 3} & \multicolumn{3}{c|}{\bf 2} & \multicolumn{1}{c|}{\bf 1} \\ \hline
1 & \gcl & \gcl & \gcl & \gcl & \gcl & \gcl  & \gcl & \gcl & \gcl & \gcl & \gcl & \gcl & \gcl & \gcl & \gcl & \gcl & \gcl & \gcl & \gcl \\
2 & \gcl & \gcl & \gcl & \gcl & \gcl & \gcl  & \gcl & \gcl & \gcl & \gcl & \gcl & \gcl & \gcl & \gcl & \gcl & \gcl & \gcl & \gcl & \gcl \\
3 & \gcl & \gcl & 7.00 & \gcl & \gcl & 11.00 & \gcl & \gcl & \gcl & \gcl & \gcl & \gcl & \gcl & \gcl & 0.00 & \gcl & \gcl & 0.00 & \gcl \\
4 & \gcl & 0.63 & 4.37 & \gcl & 2.00 & 7.00  & \gcl & 9.00 & \gcl & 6.00 & \gcl & 9.00 & \gcl & 0.00 & 0.00 & \gcl & 0.00 & 0.00 & \gcl \\
5 & 0.33 & 0.33 & 2.33 & 0.00 & 0.00 & 0.00  & 1.20 & 1.80 & 1.50 & 1.50 & 4.00 & 1.00 & 1.40 & 0.20 & 0.40 & 1.50 & 0.90 & 0.60 & 3.00 \\ \hline
\end{tabular}}
\caption{Disaggregated sales of Table \ref{ex:vvrr_op}}
\label{ex:vvrr_op_disagg}
\end{table}

After applying the EM algorithm on the disaggregated data, we end up with the disaggregated solution, presented in Table \ref{sol:vvrr_op_disagg}.

\begin{table}[!htbp]
\centering
\resizebox{1.00\linewidth}{!}{
\begin{tabular}[!htbp]{|c||cccc|cccc|c|c|ccc|ccc|cc||c|}
	\hline
\multirow{2}{*}{\bf Estimates} & \multicolumn{18}{c||}{\bf Period} & \multirow{2}{*}{$\mathbf{v}_i$} \\ \cline{2-19}
& \multicolumn{4}{c|}{\bf 15} & \multicolumn{4}{c|}{\bf 14} & \multicolumn{1}{c|}{\bf 13} & \multicolumn{1}{c|}{\bf 12} & \multicolumn{3}{c|}{\bf 11} & \multicolumn{3}{c|}{\bf 10} & \multicolumn{2}{c||}{\bf 9} & \\
\hline
1 & 0.81 & 0.91 & 1.33 & 10.00 & 2.41 & 5.02 & 2.86 & 15.00 & 11.00 & 14.00 & 6.03 & 5.27 & 15.88 & 4.05 & 2.68 & 11.50 & 0.81 & 16.85 & 1.000 \\
2 & 0.59 & 0.67 & 0.94 & 9.62 & 1.76 & 3.65 & 1.64 & 3.60 & 11.00 & 8.00 & 4.39 & 3.83 & 13.63 & 2.95 & 1.95 & 10.90 & 0.59 & 12.26 & 0.727 \\
3 & 0.24 & 0.25 & 0.38 & 3.89 & 0.71 & 1.20 & 0.91 & 2.00 & 1.00 & 11.00 & 1.77 & 1.55 & 2.73 & 1.19 & 0.90 & 2.04 & 0.24 & 6.29 & 0.294 \\
4 & 0.14 & 0.18 & 0.27 & 2.80 & 0.36 & 0.80 & 0.61 & 1.33 & 4.00 & 1.00 & 0.91 & 0.85 & 2.45 & 0.71 & 0.36 & 0.82 & 0.11 & 1.21 & 0.150 \\
5 & 0.00 & 0.00 & 0.00 & 0.00 & 0.06 & 0.36 & 0.27 & 0.60 & 0.00 & 0.00 & 0.15 & 0.07 & 0.20 & 0.00 & 0.00 & 0.00 & 0.04 & 0.40 & 0.026 \\ \hline
$\mathbf{\lambda}_t$ & 2.55 & 2.87 & 4.16 & 37.59 & 7.57 & 15.76 & 8.97 & 32.19 & 38.57 & 48.57 & 18.94 & 16.54 & 49.85 & 12.73 & 8.41 & 36.09 & 2.55 & 52.89 & \\ \hline
\end{tabular}}
\resizebox{1.00\linewidth}{!}{
\begin{tabular}[!htbp]{|c||ccc|ccc|cc|cc|cc|ccc|ccc|c||c|}
	\hline
\bf \phantom{Estimates2} & \multicolumn{3}{c|}{\bf 8} & \multicolumn{3}{c|}{\bf 7} & \multicolumn{2}{c|}{\bf 6} & \multicolumn{2}{c|}{\bf 5} & \multicolumn{2}{c|}{\bf 4} & \multicolumn{3}{c|}{\bf 3} & \multicolumn{3}{c|}{\bf 2} & \multicolumn{1}{c||}{\bf 1} & \phantom{x} $\mathbf{v}_i$ \phantom{x} \\ \hline
1 & 4.02 & 1.94 & 13.13 & 0.00 & 4.05 & 17.23 & 14.48 & 21.90 & 18.10 & 15.21 & 48.27 & 20.27 & 16.89 & 0.41 & 0.38 & 18.10 & 1.82 & 0.57 & 36.20 & 1.000 \\
2 & 2.93 & 1.41 & 9.55 & 0.00 & 2.95 & 12.54 & 10.53 & 15.93 & 13.17 & 11.06 & 35.11 & 14.75 & 12.29 & 0.29 & 0.28 & 13.17 & 1.33 & 0.42 & 26.33 & 0.727 \\
3 & 1.18 & 0.57 & 3.15 & 0.00 & 1.19 & 4.95 & 4.26 & 6.44 & 5.32 & 4.47 & 14.19 & 5.96 & 4.97 & 0.12 & 0.00 & 5.32 & 0.54 & 0.00 & 10.64 & 0.294 \\
4 & 0.60 & 0.22 & 1.97 & 0.00 & 0.71 & 3.15 & 2.17 & 3.20 & 2.72 & 2.14 & 7.24 & 3.20 & 2.53 & 0.00 & 0.00 & 2.72 & 0.00 & 0.00 & 5.43 & 0.150 \\
5 & 0.10 & 0.12 & 1.05 & 0.00 & 0.00 & 0.00 & 0.37 & 0.64 & 0.46 & 0.53 & 1.23 & 0.36 & 0.43 & 0.07 & 0.18 & 0.46 & 0.32 & 0.27 & 0.92 & 0.026 \\ \hline
$\mathbf{\lambda}_t$ & 12.62 & 6.10 & 41.19 & 0.00 & 12.73 & 54.09 & 45.45 & 68.72 & 56.81 & 47.72 & 151.49 & 63.63 & 53.02 & 1.27 & 1.20 & 56.81 & 5.73 & 1.80 & 113.62 & \\ \hline
\end{tabular}}
\caption{Estimated demand and parameters of disaggregated sales data}
\label{sol:vvrr_op_disagg}
\end{table}

Note that the estimated primary demands do not preserve the split proportions of the sales. For instance, in period 9, $\lambda_{9_1} = 2.55$ and $\lambda_{9_2} = 52.89$, while the sales were split by proportions $\alpha_1=0.1$ and $\alpha_2=0.9$. Note that the split proportions of sales do not need to match the proportion of demands, because the varying assortments imply varying market shares and spilled demand. However, to reduce the number of parameters to be estimated, it would be possible to modify the EM algorithm and incorporate constraints on $\lambda_{t_i}$ to preserve the split proportions. Aggregating the results back, we end up with the final results presented in Table \ref{sol:vvrr_op_agg}.

\begin{table}[!htbp]
\centering
\resizebox{0.95\linewidth}{!}{
\begin{tabular}[!htbp]{|c||c|c|c|c|c|c|c|c|c|c|c|c|c|c|c||c|}
	\hline
\multirow{2}{*}{\bf Estimates} & \multicolumn{15}{c||}{\bf Period} & \multirow{2}{*}{$\mathbf{v}_i$} \\ \cline{2-16}
            & \bf 15 & \bf 14 & \bf 13 & \bf 12 & \bf 11 & \bf 10 & \bf 9 & \bf 8 & \bf 7 & \bf 6 & \bf 5 & \bf 4 & \bf 3 & \bf 2 & \bf 1 & \\ \hline
1 & 13.05 & 25.29 & 11.00 & 14.00 & 27.19 & 18.23 & 17.66 & 19.09 & 21.29 & 36.38 & 33.31 & 68.54 & 17.68 & 20.50 & 36.20 & 1.000 \\
2 & 11.82 & 10.64 & 11.00 & 8.00 & 21.85 & 15.80 & 12.85 & 13.89 & 15.48 & 26.46 & 24.22 & 49.85 & 12.86 & 14.91 & 26.33 & 0.727 \\
3 & 4.76 & 4.82 & 1.00 & 11.00 & 6.05 & 4.14 & 6.53 & 4.90 & 6.14 & 10.70 & 9.79 & 20.15 & 5.09 & 5.86 & 10.64 & 0.294 \\
4 & 3.39 & 3.10 & 4.00 & 1.00 & 4.21 & 1.89 & 1.32 & 2.79 & 3.86 & 5.38 & 4.85 & 10.45 & 2.53 & 2.72 & 5.43 & 0.150 \\
5 & 0.00 & 1.29 & 0.00 & 0.00 & 0.43 & 0.00 & 0.44 & 1.27 & 0.00 & 1.01 & 1.00 & 1.59 & 0.68 & 1.05 & 0.92 & 0.026 \\  \hline
$\mathbf{\lambda}_t$ & 47.17 & 64.50 & 38.57 & 48.57 & 85.33 & 57.23 & 55.43 & 59.92 & 66.82 & 114.17 & 104.53 & 215.12 & 55.50 & 64.34 & 113.62 & \\ \hline
\end{tabular}
}
\caption{Estimated demand and parameters using data disaggregation}
\label{sol:vvrr_op_agg}
\end{table}

The results are fairly similar to the results in Table \ref{vvrr_op_sol}, where we incorporated partial availability into the attraction model. The total primary demand is 1190.8 as opposed to 1194.6. To benchmark these solutions, a naive approach would be to use the projection method to preprocess the observed sales using the open percentage information, such as
\ba
\widehat{\textnormal{sales}} = \frac{\textnormal{sales}}{\textnormal{openpct}},
\ea
and then apply EM algorithm on the preprocessed data. With this approach we estimate the overall primary demand as 1519.9. This shows that incorporating partial availability natively into the attraction model or disaggregating the sales data are much more robust approaches. In case we observe small outliers in partial availability data, we end up with very large preprocessed sales using the projection method. Experiments showed that the disaggregation method dampens the effect of outliers the most.

\subsection{Strong sell-down}
\label{sec:selldown}

The demand model discussed in \citet{vvrr12} and \citet{tarek2016} assumes that customers choose from the available products based on the Basic Attraction Model (BAM) with the IIA property. The probability of product selection is governed by (\ref{eq:BAM}). This model cannot fit to scenarios with strong sell-down, and the same applies to the Generalized Attraction Model (GAM) discussed in \citet{grs15}. Strong sell-down can be a common phenomena in practice, because customers often seek to purchase the cheapest available product. A 100\% sell-down example is presented in Table \ref{table:selldown}.

\begin{table}[!htbp]
\centering
\resizebox{0.35\linewidth}{!}{
\begin{tabular}[h]{|c||c|c|c|c|c|c|}
	\hline
\multirow{2}{*}{Sales} & \multicolumn{6}{c|}{Period} \\ \cline{2-7}
            & 6 & 5 & 4 & 3 & 2 & 1 \\
\hline
1  & 2    & 6    & 0    & 0    & 0  & 0  \\
2  & \gcl & \gcl & 13   & 15   & 0  & 0  \\
3  & \gcl & \gcl & \gcl & \gcl & 20 & 22 \\ \hline
\end{tabular}
}
\caption{Example with 100\% sell-down}
\label{table:selldown}
\end{table}

The EM algorithm developed in \citet{vvrr12} does not converge on this example, because the attraction model is unable to fit to the scenario presented in the data. The estimated first choice demands will converge to an unbounded solution.

A natural way to handle sell-down is by introducing an additional parameter $l$ for the lowest available product. The attraction model (\ref{eq:BAM}) would be extended as
\be
P_{i}(S_t, \bm v, L_t, l) = \frac{v_i + \mathbbm{1}(i \in L_t) \cdot l}{v_0 + \sum_{j \in S_t} \left(v_j + \mathbbm{1}(j \in L_t) \cdot l \right)}
\label{eq:LOC}
\ee
Parameter $l$ is an additional preference weight for the product with the lowest price, or the product with excess attraction in the assortment. In practice we could add additional preference weights by group of products, depending on the structure of the problem. $L$ is a set of the indices for the lowest available product over time, and $\mathbbm{1}(i \in L_t)$ is an indicator function, representing whether product $i$ is the lowest available product at time $t$.

For the example in Table \ref{table:selldown}, $L = (1,1,2,2,3,3)$. $L_3 = 2$ means that at $t = 3$ the lowest available class is the second class from the top.

To present a  motivating example with strong sell-down: assuming that the true parameters are $v_0 = 1$, $v_1 = 0.4$, $v_2 = 0.7$, $v_3 = 0.1$ and $l = 10$, the purchase probabilities for various assortments, induced by (\ref{eq:LOC}), are shown in Table \ref{table:loc}.

\begin{table}[!htbp]
\centering
\resizebox{0.4\linewidth}{!}{
\begin{tabular}[!htbp]{|c||c|c|c|}
	\hline
\multirow{2}{*}{Probability} & \multicolumn{3}{c|}{Assortment} \\ \cline{2-4}
            & \bf $L_1 = 1$ & \bf $L_2 = 2$ & \bf $L_3 = 3$ \\
\hline
$P_{0}$  & 8.8\%  & 8.3\%  & 8.2\%  \\
$P_{1}$  & 91.2\% & 3.3\%  & 3.3\%  \\
$P_{2}$  & \gcl   & 88.4\% & 5.7\%  \\
$P_{3}$  & \gcl   & \gcl   & 82.8\% \\ \hline
\end{tabular}
}
\caption{Purchase probabilities with additional preference for lowest available product}
\label{table:loc}
\end{table}

\noindent These purchase probabilities can describe a realistic scenario, where the lowest available product receives majority of the demand. The products available above, however, still maintain the IIA property. This extended model fits data better with strong sell-down, which cannot be described by the basic MNL model.

Essentially the same idea, with more parameters than just the scalar $l$ and a slightly more general formulation, the spiked-MNL model was considered in \citet{cao2019}. The spiked-MNL model extends the classical MNL model by having a separate attractiveness parameter for the cheapest available fare class on each flight.

The likelihood function for the extended model with parameter $l$ becomes
\footnotesize
\ba
L_I(\bm \lambda, \bm v, l) &=& \prod_{t=1}^{T} \left[ \left( P(m_t~\textnormal{customers buy in period}~t|\bm v, \bm \lambda, l)
\frac{m_t!}{z_{1t}!z_{2t}!\cdots z_{nt}!} \right) \right. \\
&& \left. \prod_{j \in S_t} \left[ \frac{P_{j}(S_t, \bm v, L_t, l)}{\sum_{i \in S_t} P_{i}(S_t, \bm v, L_t, l)} \right]^{z_{jt}} \right] \\
\ea
\normalsize
where
\footnotesize
\ba
P(m_t~\textnormal{customers buy in period}~t|\bm v, \bm \lambda, l) &=& \frac{\left[ \lambda_t \sum_{i \in S_t} P_{i}(S_t, \bm v, L_t, l) \right]^{m_t}
\exp\left(- \lambda_t \sum_{i \in S_t} P_{i}(S_t, \bm v, L_t, l) \right)}{m_t!}
\ea
\normalsize
and the log-likelihood function modifies to
\footnotesize
\be
\label{eq:loglik_loc}
l_I(\bm{v}, \bm{\lambda}, l) &=& \sum_{t=1}^{T} \left[ m_t \log\left( \frac{\lambda_t}{v_{0} + \sum_{i \in S_t} \left(v_i + \mathbbm{1}(i \in L_t) \cdot l \right)} \right) - \lambda_t \frac{\sum_{i \in S_t} \left(v_i + \mathbbm{1}(i \in L_t) \cdot l \right)}{v_{0} + \sum_{i \in S_t} \left(v_i + \mathbbm{1}(i \in L_t) \cdot l \right)} +  \right. \nonumber \\
&& \left. \sum_{i \in S_t} z_{it}\log\left(v_i + \mathbbm{1}(i \in L_t) \cdot l \right) \right]
\ee
\normalsize
Note that in this paper we are not providing an algorithm for maximizing the log-likelihood (\ref{eq:loglik_loc}), but we can rely on available nonlinear solvers, if needed. It would be of practical interest, however, to extend the EM \citep{vvrr12} and MM \citep{tarek2016} algorithms to estimate the parameters of this extended model.

\subsection{Constrained parameter space}
\label{sec:constrained_parameters}

Algorithms based on the model assumptions above can lead to unreasonable estimates in practice. For instance, on airline data we observed estimated demands which were high in a business scenario. This can happen due to data sparsity, outliers, data preprocessing steps, and most likely that the assumptions of the MNL model do not fit the true data generating model. Therefore, in practical applications, we might want to constrain the parameters of the model. \citet{tarek2016} discusses regularization as an elegant solution to the problem, and they develop algorithms for $L_1$ regularization (Lasso regression) and $L_2$ regularization (Ridge regression). It can be also of interest to apply constraints on the arrival rate of customers ($\lambda_t$), having an interpretable business fence on these parameters.

In practice we could solve a constrained maximum likelihood problem, such as
\ba
\max_{\bm \lambda, \bm v} && l_I  \\
\textrm{s.t.}&&\\
&& \sum_{t=1}^{T} \lambda_t \leq L
\ea
where the constraint ensures that the overall arrival rate or mean total demand over the time horizon is less than an upper bound $L$. $L$ could be defined, for instance, as $C$ times the total observed sales, that is
\ba
L = C \frac{1}{s} \sum_{t=1}^{T} \sum_{i=1}^{n} z_{it}
\ea
$C$ is a regularization parameter which has to be set based on practical considerations. Note that $\lambda_t$ includes the no purchase alternative, so we scale the total sales with market share $s$. Alternatively, we could constrain the arrival rates at each time period and solve constrained maximum likelihood problem
\ba
\max_{\bm \lambda, \bm v} && l_I  \\
\textrm{s.t.}&&\\
&& \lambda_t \leq L_t, \quad t = 1,\ldots,T
\ea

Similarly, we could set $L_t = C \frac{1}{s} \sum_{i=1}^{n} z_{it}$, a constant times the observed sales at time period $t$. It would be of practical interest to extend the EM algorithm of \citet{vvrr12} to solve this constrained problem. In Section \ref{sec:constrained_opt} we will extend the MM algorithm in \citet{tarek2016} to this constrained problem and also discuss how to solve the optimization problem using the Frank-Wolfe algorithm.

\subsection{Non-homogeneous product set}
\label{sec:nonhomo_prod}

In revenue management, we often encounter changing products due to changes in the system, hence we need to deal with a non-homogeneous product set over time. Instead of dividing the observed data to homogeneous parts, we want to extend the algorithms to handle non-homogeneous product offerings, allowing us to use the whole data set and borrow power to accurately estimate the parameters.

Let us look at a hypothetical airline sales example in Table \ref{table:sch_ex1}, which was created from the synthetic example in \citet{vvrr12}. We did this, so the estimated results are easy to compare.

\begin{table}[!htbp]
\centering
\resizebox{1.00\linewidth}{!}{
\begin{tabular}[h]{|cc||c|c|c|c|c|c|c|c|c|c|c|c|c|c|c|c|c|c|c|c|c|c|c|c|c|c|c|c|c|c|}
	\hline
\multicolumn{2}{|c||}{\multirow{2}{*}{\bf Sales}} & \multicolumn{30}{c|}{\bf Period} \\ \cline{3-32}
            & & \bf 1 & \bf 2 & \bf 3 & \bf 4 & \bf 5 & \bf 6 & \bf 7 & \bf 8 & \bf 9 & \bf 10 & \bf 11 & \bf 12 & \bf 13 & \bf 14 & \bf 15 & \bf 16 & \bf 17 & \bf 18 & \bf 19 & \bf 20 & \bf 21 & \bf 22 & \bf 23 & \bf 24 & \bf 25 & \bf 26 & \bf 27 & \bf 28 & \bf 29 & \bf 30 \\
\hline
\textnormal{1} & \textnormal{flt1-prod1} & 10 & 15 & 11 & 14 &\dcl&\dcl&\dcl&\dcl&\dcl&\dcl&\dcl&\dcl&\dcl&\dcl&\dcl& \multicolumn{15}{c|}{\gcl} \\
\textnormal{2} & \textnormal{flt1-prod2} & 11 & 6  & 11 & 8  & 20 & 16 &\dcl&\dcl&\dcl&\dcl&\dcl&\dcl&\dcl&\dcl&\dcl& \multicolumn{15}{c|}{\gcl} \\
\textnormal{3} & \textnormal{flt1-prod3} & 5  & 6  & 1  & 11 & 4  & 5  & 14 & 7  & 11 &\dcl&\dcl&\dcl&\dcl&\dcl&\dcl& \multicolumn{15}{c|}{\gcl} \\
\textnormal{4} & \textnormal{flt1-prod4} & 4  & 4  & 4  & 1  & 6  & 4  & 3  & 5  & 9  & 9  & 6  & 9  &\dcl&\dcl&\dcl& \multicolumn{15}{c|}{\gcl} \\
\textnormal{5} & \textnormal{flt1-prod5} & 0  & 2  & 0  & 0  & 1  & 0  & 1  & 3  & 0  & 3  & 3  & 5  & 2  & 3  & 3  & \multicolumn{15}{c|}{\gcl} \\ \hline
\textnormal{6} & \textnormal{flt2-prod1} & \multicolumn{15}{c|}{\gcl}  & 10 & 15 & 11 & 14 &\dcl&\dcl&\dcl&\dcl&\dcl&\dcl&\dcl&\dcl&\dcl&\dcl&\dcl \\
\textnormal{7} & \textnormal{flt2-prod2} & \multicolumn{15}{c|}{\gcl}  & 11 & 6  & 11 & 8  & 20 & 16 &\dcl&\dcl&\dcl&\dcl&\dcl&\dcl&\dcl&\dcl&\dcl \\
\textnormal{8} & \textnormal{flt2-prod3} & \multicolumn{15}{c|}{\gcl}  & 5  & 6  & 1  & 11 & 4  & 5  & 14 & 7  & 11 &\dcl&\dcl&\dcl&\dcl&\dcl&\dcl \\
\textnormal{9} & \textnormal{flt2-prod4} & \multicolumn{15}{c|}{\gcl}  & 4  & 4  & 4  & 1  & 6  & 4  & 3  & 5  & 9  & 9  & 6  & 9  &\dcl&\dcl&\dcl \\
\textnormal{10} & \textnormal{flt2-prod5} & \multicolumn{15}{c|}{\gcl}  & 0  & 2  & 0  & 0  & 1  & 0  & 1  & 3  & 0  & 3  & 3  & 5  & 2  & 3  & 3   \\ \hline
\textnormal{11} & \textnormal{flt3-prod1} & 20 & 30 & 22 & 28 &\dcl&\dcl&\dcl&\dcl&\dcl&\dcl&\dcl&\dcl&\dcl&\dcl&\dcl& 20 & 30 & 22 & 28 &\dcl&\dcl&\dcl&\dcl&\dcl&\dcl&\dcl&\dcl&\dcl&\dcl&\dcl \\
\textnormal{12} & \textnormal{flt3-prod2} & 22 & 12 & 22 & 16 & 40 & 32 &\dcl&\dcl&\dcl&\dcl&\dcl&\dcl&\dcl&\dcl&\dcl& 22 & 12 & 22 & 16 & 40 & 32 &\dcl&\dcl&\dcl&\dcl&\dcl&\dcl&\dcl&\dcl&\dcl \\
\textnormal{13} & \textnormal{flt3-prod3} & 10 & 12 & 2  & 22 & 8  & 10 & 28 & 14 & 22 &\dcl&\dcl&\dcl&\dcl&\dcl&\dcl& 10 & 12 & 2  & 22 & 8  & 10 & 28 & 14 & 22 &\dcl&\dcl&\dcl&\dcl&\dcl&\dcl \\
\textnormal{14} & \textnormal{flt3-prod4} & 8  & 8  & 8  & 2  & 12 & 8  & 6  & 10 & 18 & 18 & 12 & 18 &\dcl&\dcl&\dcl& 8  & 8  & 8  & 2  & 12 & 8  & 6  & 10 & 18 & 18 & 12 & 18 &\dcl&\dcl&\dcl \\
\textnormal{15} & \textnormal{flt3-prod5} & 0  & 4  & 0  & 0  & 2  & 0  & 2  & 6  & 0  & 6  & 6  & 10 & 4  & 6  & 6  & 0  & 4  & 0  & 0  & 2  & 0  & 2  & 6  & 0  & 6  & 6  & 10 & 4  & 6  & 6  \\ \hline
\end{tabular}
}
\caption{Schedule change example}
\label{table:sch_ex1}
\end{table}

We have 3 flights, each of them having 5 different products, a total of $n=15$ products. The products of flight 1 only exist at the first 15 time periods, and after that the label changes to flight 2. This can happen for various reasons, for example schedule change. It is not always possible to preprocess the data and match the products of flight 1 to flight 2. We can see that the product offerings are non-homogeneous over time, or unbalanced. Product 1 of flight 1 is not available for sale for time periods 5-15, which is distinguished from not being part of the product offer set at time periods 16-30. Notice that the observed sales data of flight 1 for time periods 1-15 is the same as the synthetic example created in \citet{vvrr12}, and the observed sales for Flight 3 are just twice of that, like a flight with double the demand and capacity.

To introduce non-homogeneous product set into the notation, let $I_t$ denote the set of existing products of the retailer at time $t$. Note that this is different from $S_t$, which is the set of products available for purchase at time $t$, therefore $S_t \subseteq I_t$. For instance, in Table \ref{table:sch_ex1}, $I_7 = \{1,2,3,4,5,11,12,13,14,15\}$ and $S_7 = \{3,4,5,13,14,15\}$. The set of all products is $I = \bigcup_{t=1}^T I_t = \left\{1,2,\ldots,15\right\}$.

In Section \ref{sec:extended_em} we will extend the EM algorithm of \citet{vvrr12} for non-homogeneous product offerings. The MM and Frank-Wolfe algorithms discussed in Section \ref{sec:constrained_opt} will also be able to natively handle this product structure.

\subsection{Handling market share constraint and no-purchase option}
\label{sec:share_nopurchase}

We mentioned that the objective function (\ref{eq:loglik}) has an infinite number of maximum likelihood estimates, therefore an additional constraint is applied in \citet{vvrr12} on the preference weights as a function of the market share
\be
\label{eq:share_constr}
s = \frac{\sum_{i=1}^n v_i}{v_0 + \sum_{i=1}^n v_i}
\ee
using the standard scaling of $v_0 = 1$. The market share constraint is linear but the objective function is non-convex, therefore \citet{tarek2016} formulates the problem using $v_i = \exp(\beta_i)$ in the objective function, so the market share constraint becomes
\ba
s = \frac{\sum_{i=1}^n \exp(\beta_i)}{1 + \sum_{i=1}^n \exp(\beta_i)}.
\ea
They establish that the objective function in this space is a concave function, and they show how to deal with the non-convex market share constraint by solving the problem without the constraint and then using a transformation to satisfy the constraint.

The outside alternative could represent the competitor's product, or both the competitor's product and the no-purchase option, which would lead to different interpretations of $s$ and $\bm \lambda$ \citep{vvrr12}. The outside alternative is treated as a separate product that is always available, therefore $v_0$ is constant and the standard scaling $v_0 = 1$ is used.

To look further, let $I$ denote the set of all products of the retailer, that is $I = \left\{1,2,\ldots,n\right\}$, so it follows that $S_t \subseteq I$. With the market share constraint above, the probability of selecting one of the retailer's product, when everything is available for sale, is $s$, that is
\ba
P_I(I, \bm v) = \frac{\sum_{j \in I} v_j}{v_0 + \sum_{j \in I} v_j} = s
\ea
When the retailer offers an assortment $S_t$ with a subset of all the products, it follows that
\ba
\tilde{s}_t = P_{S_t}(S_t, \bm v) = \frac{\sum_{j \in S_t} v_j}{v_0 + \sum_{j \in S_t} v_j} < s
\ea
given $v_0$ is fixed. This means that the model induced market share ($\tilde{s}_t$) of the retailer at time $t$ is less than $s$, and the share of outside alternative is larger than $1-s$. If we think of the outside alternative as the competitor's product, one could argue that the competitor's product might not always be available for purchase either, so the preference weight $v_0$ is not constant over time. Also, the estimated market share $s$ is calculated over all possible sets of assortments, not only when all the retailer's products are available for sale.

To extend this model, we can relax the assumption of constant, time-independent outside alternative and introduce $v_{0t}$, allowing the preference weight to change over time. We also need this extension to be able to incorporate non-homogeneous product sets into the market share constraint, by using $I_t$, the offer set of retailer at time $t$ (Section \ref{sec:nonhomo_prod}). The market share constraint modifies to a set of constraints
\be
\label{eq:share_constr1}
s = \frac{\sum_{i \in I_t} v_i}{v_{0t} + \sum_{i \in I_t} v_i},~~t=1,\ldots,n
\ee
which ensure that the share of outside alternative is $1-s$, for all time points $t$, when all the retailer's products are available for sale. That is
\ba
P_{0}(I, \bm v) = \frac{v_{0t}}{v_{0t} + \sum_{i \in I_t} v_i} = \frac{\frac{1-s}{s}\sum_{i \in I_t} v_i}{\frac{1-s}{s}\sum_{i \in I_t} v_i + \sum_{i \in I_t} v_i} = 1-s.
\ea

Now let us consider an edge case scenario where the retailer's model induced market share is constant over time regardless of the assortment he offers. This could be achieved with the set of market share constraints
\be
\label{eq:share_constr2}
s = \frac{\sum_{i \in S_t} v_i}{v_{0t} + \sum_{i \in S_t} v_i},~~t=1,\ldots,n
\ee
leading to
\ba
\tilde{s}_t = P_{S_t}(S_t, \bm v) = \frac{\sum_{j \in S_t} v_j}{v_{0t} + \sum_{j \in S_t} v_j} = \frac{\sum_{j \in S_t} v_j}{\frac{1-s}{s}\sum_{j \in S_t} v_j + \sum_{j \in S_t} v_j} = s.
\ea
The retailers's share remains $s$ and the share of the outside alternative remains $1-s$, at each time $t$, regardless of what products are available for purchase. This edge case might not make sense for all applications, but it is interesting to consider as an alternative to the assumption that the outside alternative is always fully available. In the airline retail context we could think of the outside alternative as the competitor's products whose availability can be as limited as the host airline's products.

To introduce a continuum of cases between (\ref{eq:share_constr1}) and (\ref{eq:share_constr2}), let us introduce parameter $\alpha \in [0,1]$ controlling the availability of outside alternative. $\alpha = 0$ represents the case where the outside alternative is always available for sale (\ref{eq:share_constr1}), and naturally, $\alpha=1$ represents the case where the outside alternative limits availability the same fashion as the retailer. This is more of a hypothetical case, given the outside alternative could include the no purchase option, which is always an available choice. Combining (\ref{eq:share_constr1}) and (\ref{eq:share_constr2}), the set of market share constraints become
\be
\label{eq:share_constr_alpha}
s = (1-\alpha) \frac{\sum_{i \in I_t} v_i}{v_{0t} + \sum_{i \in I_t} v_i} + \alpha \frac{\sum_{i \in S_t} v_i}{v_{0t} + \sum_{i \in S_t} v_i},~~t=1,\ldots,n
\ee
Using $\alpha=0$, $I_t = I = \{1,\ldots,n\}$ and $v_{0t}=1$, the constraints simplify to (\ref{eq:share_constr}).

The models require the knowledge of market share $s$, which in practice can be difficult to acquire, and the estimate itself can be inaccurate. \citet{tarek2016} included a study on the sensitivity of model estimates to the input market share. We mentioned that, in practice, market share $s$ is most likely estimated from data observed over various sets of assortments, not only when all the retailer's products are available for sale. Therefore, it could also make sense to use a market share constraint aggregated over time horizon $T$, that is
\be
\label{eq:share_constr_agg}
s = \frac{\sum_{t=1}^T \sum_{i \in I_t} v_i}{ \sum_{t=1}^T \left( v_{0t} + \sum_{i \in I_t} v_i \right)}
\ee
This constraint ensures that the market share over the time horizon $T$ is $s$, taking into account the changing offer set $I_t$ over time. This is less restrictive and could be a more reasonable assumption than using (\ref{eq:share_constr}), which forces the market share at each time point $t$ with a set of constraints. Adding $\alpha$ to control availability of the outside alternative, we would use
\be
\label{eq:share_constr_agg_alpha}
s = \frac{\sum_{t=1}^T \left[ (1-\alpha) \sum_{i \in I_t} v_i + \alpha \sum_{i \in S_t} v_i \right]}{ \sum_{t=1}^T \left[ v_{0t} + (1-\alpha) \sum_{i \in I_t} v_i + \alpha \sum_{i \in S_t} v_i \right]}
\ee

Finally, we want to point out that instead of solving the market share constrained optimization problem
\ba
\max_{\bm \lambda, \bm v} && \sum_{t=1}^{T} \left[ m_t \log\left( \frac{\lambda_t}{v_{0} + \sum_{i \in S_t} v_i} \right) - \lambda_t \frac{\sum_{i \in S_t} v_i}{v_{0} + \sum_{i \in S_t} v_i} +\sum_{i \in S_t} z_{it}\log(v_i) \right] \\
\textrm{s.t.}&&\\
&& s = \frac{\sum_{i=1}^n v_i}{v_0 + \sum_{i=1}^n v_i} \\ \\
&& v_0 = 1
\ea
we can directly incorporate the market share constraint into the objective function through $v_0 = r \sum_{i = 1}^n v_i$, with $r = (1-s)/s$, and solve
\ba
\max_{\bm \lambda, \bm v} && \sum_{t=1}^{T} \left[ m_t \log\left( \frac{\lambda_t}{r \sum_{i = 1}^n v_i + \sum_{i \in S_t} v_i} \right) - \lambda_t \frac{\sum_{i \in S_t} v_i}{r \sum_{i = 1}^n v_i + \sum_{i \in S_t} v_i} +\sum_{i \in S_t} z_{it}\log(v_i) \right] \\
\textrm{s.t.}&&\\
&& \sum_{i=1}^n v_i = 1
\ea
or by setting $v_1 = 1$. The scaling constraint on $\bm{v}$ is required to avoid multiple solutions.

Incorporating non-homogeneous product set and control of availability of outside alternative, we would need to solve
\ba
\max_{\bm \lambda, \bm v} && \sum_{t=1}^{T} \left[ m_t \log\left( \frac{\lambda_t}{v_{0t} + \sum_{i \in S_t} v_i} \right) - \lambda_t \frac{\sum_{i \in S_t} v_i}{v_{0t} + \sum_{i \in S_t} v_i} +\sum_{i \in S_t} z_{it}\log(v_i) \right] \\
\textrm{s.t.}&&\\
&& v_{0t} = r \left[ (1-\alpha) \sum_{i \in I_t} v_i + \alpha \sum_{i \in S_t} v_i \right] \\ \\
&& \sum_{i=1}^n v_i = 1
\ea
We will see in Section \ref{sec:frank_wolfe} that this formulation with upper bound constraints on the arrival rate is easy to solve using the Frank-Wolfe method.

\section{Extended EM algorithm}
\label{sec:extended_em}

In this section we extend the EM algorithm of \citet{vvrr12} with non-homogeneous product sets (Section \ref{sec:nonhomo_prod}) and the ability to control the availability of outside alternative (Section \ref{sec:share_nopurchase}).

We incorporate $I_t$, the offer set of retailer at time $t$, and $v_{0t}$, the time-dependent preference weight for the outside alternative into the algorithm. We use
\ba
v_{0t} = r \left[ (1-\alpha) \sum_{i \in I_t} v_i + \alpha \sum_{i \in S_t} v_i \right]
\ea
where $r = \frac{1-s}{s}$ and $\alpha \in [0,1]$ controls the availability of the outside alternative. Adding these into the formulation, we can modify the key E-step equations of the EM algorithm as
\ba
\hat{X}_{jt} &=&
\begin{cases}
    \frac{v_j}{\left(1+r\right)\sum_{i \in I_t} v_i} \frac{\sum_{h \in S_t} v_h + v_{0t}}{\sum_{h \in S_t} v_h}
             \sum_{h \in S_t} z_{ht}, & \text{if $j \notin S_t \cup \{0\}$} \\
    \frac{\sum_{h \in S_t} v_h + v_{0t}} {\left(1+r\right) \sum_{i \in I_t} v_i} z_{jt}, & \text{if $j \in S_t$}
\end{cases} \\
\hat{Y}_{jt} &=& \frac{\sum_{h \notin S_t \cup \{0\}} v_h - v_{0t} + r \sum_{i \in I_t} v_i} {\left(1+r\right) \sum_{i \in I_t} v_i} z_{jt},
~j \in S_t \\
\hat{X}_{0t} &=& r \sum_{i \in I_t} \hat{X}_{it} \\
\hat{Y}_{0t} &=& \frac{v_{0t}}{\sum_{i \in S_t} v_i + v_{0t}}  \sum_{h \notin S_t \cup \{0\}} \hat{X}_{ht}
\ea

For the M-step we need to maximize the conditional expected, complete data log-likelihood function, which becomes
\be
{\mathcal{L}}(\bm v) = \sum_{t = 1}^T \sum_{j \in I_t} \hat{X}_{jt} \log \left( \frac{v_j}{(1+r)\sum_{i \in I_t} v_i} \right) + \sum_{t = 1}^T \hat{X}_{0t} \log \left( s \right)\label{eq:compLL}
\ee
To find the maxima, the first order conditions become
\ba
\frac{\partial {\mathcal{L}}}{\partial v_i} = \sum_{t=1}^{T} \mathbbm{1} \left[i \in I_t \right] \left( \frac{\hat{X}_{it}}{v_i} - \frac{\sum_{k \in I_t} \hat{X}_{kt}}{\sum_{k \in I_t} v_k} \right) = 0, ~~i=1,\ldots,n
\ea
which is a system of $n$ nonlinear equations in $v_i$. The solution can be obtained with existing implementations of iterative methods, such as Newton-Raphson. We can also use a simple fixed point iteration algorithm by rewriting $\nabla \mathcal{L} = F(\bm v) = 0$ to $\Phi\left(\bm v\right) = \bm v$ and then using the fixed point iteration $\bm v^{(l+1)} := \Phi\left(\bm v^{(l)}\right)$. In our case the update becomes
\ba
v_i^{(l+1)} := \frac{\sum_{t=1}^{T} \mathbbm{1} \left[i \in I_t \right] X_{it}}{\sum_{t=1}^{T} \mathbbm{1} \left[i \in I_t \right] \frac{\sum_{k \in I_t} X_{kt}}{\sum_{k \in I_t} v_k^{(l)}}}, ~~i=1,\ldots,n
\ea

In case of homogeneous product set, so that $I_t = \{1,\ldots,n\},~\forall t$, we get a closed form solution to the M-step, that is
\ba
v_i = \frac{\sum_{t=1}^{T} X_{it}}{\sum_{j=1}^{n} \sum_{t=1}^{T} X_{jt}}, ~~ i=1,\ldots,n
\ea
The result is easy to obtain by assuming $\sum_{i=1}^{n} v_i = 1$. Note that this closed-form equation is different from the one derived in \citet{vvrr12}, since they used a different parametrization by setting $v_0=1$.

Algorithm \ref{algo_mfem_ext} presents the extended EM algorithm with non-homogeneous product set and the ability to control the availability of outside alternative.

\begin{algorithm}[!htbp]
\footnotesize 
\caption{EM algorithm with non-homogeneous product set and the control of availability of outside alternative (OA)}
\label{algo_mfem_ext}
\begin{algorithmic}[1]
\State $I_t$: set of offered products at time $t$ (product set)
\State $S_t$: set of products available for purchase at time $t$ ($S_t \subset I_t$)
\State $r = \frac{1-s}{s}$, where $s$ is market share of the retailer
\State $z_{jt}$: observed sales at time $t$ for product $j$
\State $\alpha$: control parameter ($\alpha = 1$: OA available as retailer; $\alpha = 0$: OA fully available)
\State \textbf{E-step}
\For{$t = 1,\ldots,T$}
    \State $v_{0t} = r \left[ \left(1-\alpha\right) \sum_{i \in I_t} v_i + \alpha \sum_{h \in S_t} v_h \right]$
    \Comment OA preference weight at $t$
    \For{$j \in I_t$}
        \If{$j \notin S_t$}
            \State $X_{jt} = \frac{v_j}{\left(1+r\right)\sum_{i \in I_t} v_i} \frac{\sum_{h \in S_t} v_h + v_{0t}}{\sum_{h \in S_t} v_h}
             \sum_{h \in S_t} z_{ht} $
        \Else{ ($j \in S_t$)}
        \State $Y_{jt} = \frac{\sum_{h \notin S_t} v_h - v_{0t} + r \sum_{i \in I_t} v_i} {\left(1+r\right) \sum_{i \in I_t} v_i} z_{jt}$
        \State $X_{jt} = z_{jt} - Y_{jt}$
        \EndIf
    \EndFor
    \State $X_{ot} = r \sum_{i \in I_t} X_{it}$
    \State $Y_{ot} = \frac{v_{0t}}{\sum_{h \in S_t} v_h + v_{0t}} \sum_{h \notin S_t} X_{ht}$
\EndFor
\State \textbf{M-step}
\State Find $v_i, ~ i=1,\ldots,n$ as solution to the system of nonlinear equations $F(\bm v) = 0$: \State $\sum_{t=1}^{T} \mathbbm{1} \left[i \in I_t \right] \left( \frac{X_{it}}{v_i} - \frac{\sum_{k \in I_t} X_{kt}}{\sum_{k \in I_t} v_k} \right) = 0, ~~ i=1,\ldots,n$
\State Special case (homogeneous product set), if $I_t = \{1,\ldots,n\}, ~\forall t$:
\State $v_i = \frac{\sum_{t=1}^{T} X_{it}}{\sum_{j=1}^{n} \sum_{t=1}^{T} X_{jt}}, ~~ i=1,\ldots,n$
\end{algorithmic}
\end{algorithm}

If, in practice, we observe data with partial availability, we recommend to split the sales to fully open and closed assortments (Algorithm \ref{algo_sales_split}), apply extended EM (Algorithm \ref{algo_mfem_ext}), and aggregate the solution back. It would be an interesting future research topic to further extend the EM algorithm and incorporate constraints on the arrival rates (Section \ref{sec:constrained_parameters}).

\subsection{Example: limited OA availability}

In this section we apply Algorithm \ref{algo_mfem_ext} on the simulated example from \citet{vvrr12} and demonstrate what happens when we limit the availability of the outside alternative simultaneously with the retailer's availability ($\alpha = 1$). The observed sales data is presented in Table \ref{table:vvrr_ex}.

\begin{table}[!htbp]
\centering
\resizebox{0.6\linewidth}{!}{
\begin{tabular}[h]{|c||c|c|c|c|c|c|c|c|c|c|c|c|c|c|c|}
	\hline
\multirow{2}{*}{\bf Sales} & \multicolumn{15}{c|}{\bf Period} \\ \cline{2-16}
            & \bf 15 & \bf 14 & \bf 13 & \bf 12 & \bf 11 & \bf 10 & \bf 9 & \bf 8 & \bf 7 & \bf 6 & \bf 5 & \bf 4 & \bf 3 & \bf 2 & \bf 1 \\
\hline
1  & 10 & 15 & 11 & 14 &\dcl&\dcl&\dcl&\dcl&\dcl&\dcl&\dcl&\dcl&\dcl&\dcl&\dcl \\
2  & 11 & 6  & 11 & 8  & 20 & 16 &\dcl&\dcl&\dcl&\dcl&\dcl&\dcl&\dcl&\dcl&\dcl \\
3  & 5  & 6  & 1  & 11 & 4  & 5  & 14 & 7  & 11 &\dcl&\dcl&\dcl&\dcl&\dcl&\dcl \\
4  & 4  & 4  & 4  & 1  & 6  & 4  & 3  & 5  & 9  & 9  & 6  & 9  &\dcl&\dcl&\dcl \\
5  & 0  & 2  & 0  & 0  & 1  & 0  & 1  & 3  & 0  & 3  & 3  & 5  & 2  & 3  & 3   \\ \hline
\end{tabular}
}
\caption{Simulated example of \citet{vvrr12}}
\label{table:vvrr_ex}
\end{table}

First let us look at the solution, using $\alpha = 0$, presented in Table \ref{table:vvrr_sol_alpha0}. We closely recover the results of \citet{vvrr12}, since we make the same assumption that the outside alternative is always available.

\begin{table}[!htbp]
\centering
\resizebox{0.95\linewidth}{!}{
\begin{tabular}[!htbp]{|c||c|c|c|c|c|c|c|c|c|c|c|c|c|c|c||c|}
	\hline
\multirow{2}{*}{\bf Estimates} & \multicolumn{15}{c||}{\bf Period} & \multirow{2}{*}{$\mathbf{v}_i$} \\ \cline{2-16}
            & \bf 15 & \bf 14 & \bf 13 & \bf 12 & \bf 11 & \bf 10 & \bf 9 & \bf 8 & \bf 7 & \bf 6 & \bf 5 & \bf 4 & \bf 3 & \bf 2 & \bf 1 & \\ \hline
1  & 10 & 15 & 11 & 14 & 15.03 & 12.12 & 13.04 & 10.87 & 14.49 & 15.91 & 11.93 & 18.56 & 11.51 & 17.27 & 17.27 & 1.000 \\
2  & 11 &  6 & 11 & 8  & 14.35 & 11.48 & 10.45 &  8.71 & 11.61 & 12.75 &  9.56 & 14.88 &  9.23 & 13.84 & 13.84 & 0.801 \\
3  &  5 &  6 &  1 & 11 &  2.87 &  3.59 &  6.88 &  3.44 &  5.41 &  6.22 &  4.67 &  7.26 &  4.50 &  6.75 &  6.75 & 0.391 \\
4  &  4 &  4 &  4 &  1 &  4.31 &  2.87 &  1.47 &  2.46 &  4.42 &  3.43 &  2.29 &  3.43 &  2.68 &  4.02 &  4.02 & 0.233 \\
5  &  0 &  2 &  0 &  0 &  0.72 &  0.00 &  0.49 &  1.47 &  0.00 &  1.14 &  1.14 &  1.91 &  0.63 &  0.95 &  0.95 & 0.055 \\ \hline
$\mathbf{\lambda}_t$ & 42.86 & 47.14 & 38.57 & 48.57 & 53.26 & 42.95 & 46.19 & 38.50 & 51.33 & 56.37 & 42.28 & 65.76 & 40.78 & 61.18 & 61.18 & \\ \hline
\end{tabular}
}
\caption{Estimated demand and parameters using EM algorithm with $\alpha=0$}
\label{table:vvrr_sol_alpha0}
\end{table}

Using $\alpha = 1$, the estimated primary demand at each time period is equal to the observed purchases, because the availability of outside alternative is restricted as the retailer's availability, preserving the model induced market share $s$ at each time $t$. The results are presented in Table \ref{table:vvrr_sol_alpha1}.

\begin{table}[!htbp]
\centering
\resizebox{0.95\linewidth}{!}{
\begin{tabular}[!htbp]{|c||c|c|c|c|c|c|c|c|c|c|c|c|c|c|c||c|}
	\hline
\multirow{2}{*}{\bf Estimates} & \multicolumn{15}{c||}{\bf Period} & \multirow{2}{*}{$\mathbf{v}_i$} \\ \cline{2-16}
            & \bf 15 & \bf 14 & \bf 13 & \bf 12 & \bf 11 & \bf 10 & \bf 9 & \bf 8 & \bf 7 & \bf 6 & \bf 5 & \bf 4 & \bf 3 & \bf 2 & \bf 1 & \\ \hline
1  & 10 & 15 & 11 & 14 & 12.50 & 10.08 & 7.26 & 6.05 & 8.06 & 4.84 & 3.63 & 5.65 & 0.81 & 1.21 & 1.21 & 1.000 \\
2  & 11 & 6  & 11 & 8  & 11.94 & 9.55 & 5.75 & 4.79 & 6.39 & 3.83 & 2.87 & 4.47 & 0.64 & 0.96 & 0.96 & 0.792 \\
3  & 5  & 6  & 1  & 11 & 2.39 & 2.98 & 3.88 & 1.94 & 3.05 & 1.92 & 1.44 & 2.24 & 0.32 & 0.48 & 0.48 & 0.396 \\
4  & 4  & 4  & 4  & 1  & 3.58 & 2.39 & 0.83 & 1.39 & 2.50 & 1.06 & 0.71 & 1.06 & 0.20 & 0.30 & 0.30 & 0.245 \\
5  & 0  & 2  & 0  & 0  & 0.60 & 0.00 & 0.28 & 0.83 & 0.00 & 0.35 & 0.35 & 0.59 & 0.04 & 0.06 & 0.06 & 0.046 \\ \hline
$\mathbf{\lambda}_t$ & 42.86 & 47.14 & 38.57 & 48.57 & 44.29 & 35.71 & 25.71 & 21.43 & 28.57 & 17.14 & 12.86 & 20.00 & 2.86 & 4.29 & 4.29 & \\ \hline
\end{tabular}
}
\caption{Estimated demand and parameters using EM algorithm with $\alpha=1$}
\label{table:vvrr_sol_alpha1}
\end{table}


Notice that the estimated preference weights in Tables \ref{table:vvrr_sol_alpha0} and \ref{table:vvrr_sol_alpha1} are close to each other, the estimates are not sensitive to the value of $\alpha$. However, the estimated arrival rates changed drastically, due to the change in the model induced market shares.

Parameter $\alpha$ can be used as a tool for the practitioner to control the availability of the outside alternative between these two edge cases. We would like to emphasize, however, that competitor matching can be risky and the value of $\alpha$ should ideally be inferred from exogenous data. The conservative practice is to use $\alpha = 0$.

\subsection{Example: non-homogeneous product set}
\label{sec:sch_change}

To demonstrate the extended EM algorithm on a non-homogeneous product set, let us revisit the example presented in Table \ref{table:sch_ex1} (Section \ref{sec:nonhomo_prod}). This is a hypothetical airline sales example with a schedule change. We observe 3 flights, where the sales of products of flight 1 are discontinued and the repeated as flight 2, and flight 3 has twice the observed sales of flights 1 and 2. The solution using $\alpha = 0$ is presented in Table \ref{table:sch_ex1_sol}.

\begin{table}[!htbp]
\centering
\resizebox{0.95\linewidth}{!}{
\begin{tabular}[!htbp]{|c||c|c|c|c|c|c|c|c|c|c|c|c|c|c|c||c|}
	\hline
\multirow{2}{*}{\bf Estimates} & \multicolumn{15}{c||}{\bf Period} & \multirow{2}{*}{$\mathbf{v}_i$} \\ \cline{2-16}
& \bf 1 & \bf 2 & \bf 3 & \bf 4 & \bf 5 & \bf 6 & \bf 7 & \bf 8 & \bf 9 & \bf 10 & \bf 11 & \bf 12 & \bf 13 & \bf 14 & \bf 15 & \\ \hline
\textnormal{flt1-prod1} & 10    & 15    & 11    & 14    & 15.03 & 12.12 & 13.04 & 10.87 & 14.49 & 15.91 & 11.93 & 18.56 & 11.51 & 17.27 & 17.27
                        & 1.000 \\
\textnormal{flt1-prod2} & 11    &  6    & 11    &  8    & 14.35 & 11.48 & 10.45 &  8.71 & 11.61 & 12.75 &  9.56 & 14.88 &  9.23 & 13.84 & 13.84
                        & 0.801 \\
\textnormal{flt1-prod3} & 5     &  6    &  1    & 11    &  2.87 &  3.59 &  6.88 &  3.44 &  5.41 &  6.22 &  4.67 &  7.26 &  4.50 &  6.75 &  6.75
                        & 0.391 \\
\textnormal{flt1-prod4} & 4     &  4    &  4    &  1    &  4.31 &  2.87 &  1.47 &  2.46 &  4.42 &  3.43 &  2.29 &  3.43 &  2.68 &  4.02 &  4.02
                        & 0.233 \\
\textnormal{flt1-prod5} & 0     &  2    &  0    &  0    &  0.72 &  0.00 &  0.49 &  1.47 &  0.00 &  1.14 &  1.14 &  1.91 &  0.63 &  0.95 &  0.95
                        & 0.055 \\ \hline
\textnormal{flt2-prod1} & \multicolumn{15}{c|}{\gcl}
                        & 1.000 \\
\textnormal{flt2-prod2} & \multicolumn{15}{c|}{\gcl}
                        & 0.801 \\
\textnormal{flt2-prod3} & \multicolumn{15}{c|}{\gcl}
                        & 0.391 \\
\textnormal{flt2-prod4} & \multicolumn{15}{c|}{\gcl}
                        & 0.233 \\
\textnormal{flt2-prod5} & \multicolumn{15}{c|}{\gcl}
                        & 0.055 \\ \hline
\textnormal{flt3-prod1} & 20    & 30    & 22    & 28    & 30.07 & 24.25 & 26.08 & 21.73 & 28.98 & 31.82 & 23.87 & 37.12 & 23.02 & 34.54 & 34.54
                        & 2.000 \\
\textnormal{flt3-prod2} & 22    & 12    & 22    & 16    & 28.71 & 22.97 & 20.90 & 17.42 & 23.22 & 25.50 & 19.13 & 29.75 & 18.45 & 27.68 & 27.68
                        & 1.603 \\
\textnormal{flt3-prod3} & 10    & 12    &  2    & 22    &  5.74 &  7.18 & 13.76 &  6.88 & 10.81 & 12.44 &  9.33 & 14.52 &  9.00 & 13.51 & 13.51
                        & 0.782 \\
\textnormal{flt3-prod4} &  8    &  8    &  8    &  2    &  8.61 &  5.74 &  2.95 &  4.92 &  8.85 &  6.86 &  4.57 &  6.86 &  5.36 &  8.04 &  8.04
                        &  0.465 \\
\textnormal{flt3-prod5} &  0    &  4    &  0    &  0    &  1.44 &  0.00 &  0.98 &  2.95 &  0.00 &  2.29 &  2.29 &  3.81 &  1.26 &  1.89 &  1.89
                        &  0.110 \\ \hline
$\mathbf{\lambda}_t$    & 128.57& 141.43& 115.71& 145.71& 159.79& 128.86& 138.58& 115.49& 153.98& 169.10& 126.83& 197.29& 122.35& 183.53& 183.53 &       \\ \hline
\end{tabular}}
\resizebox{0.95\linewidth}{!}{
\begin{tabular}[!htbp]{|c||c|c|c|c|c|c|c|c|c|c|c|c|c|c|c||c|}
	\hline
\bf \phantom{Estimates} & \bf 16 & \bf 17 & \bf 18 & \bf 19 & \bf 20 & \bf 21 & \bf 22 & \bf 23 & \bf 24 & \bf 25 & \bf 26 & \bf 27 & \bf 28 & \bf 29 & \bf 30 & \phantom{$\mathbf{v}_i$} \\ \hline
\textnormal{flt1-prod1} & \multicolumn{15}{c|}{\gcl}    & 1.000 \\
\textnormal{flt1-prod2} & \multicolumn{15}{c|}{\gcl}    & 0.801 \\
\textnormal{flt1-prod3} & \multicolumn{15}{c|}{\gcl}    & 0.391 \\
\textnormal{flt1-prod4} & \multicolumn{15}{c|}{\gcl}    & 0.233 \\
\textnormal{flt1-prod5} & \multicolumn{15}{c|}{\gcl}    & 0.055 \\ \hline
\textnormal{flt2-prod1} & 10    & 15    & 11    & 14    & 15.03 & 12.12 & 13.04 & 10.87 & 14.49 & 15.91 & 11.93 & 18.56 & 11.51 & 17.27 & 17.27 & 1.000 \\
\textnormal{flt2-prod2} & 11    &  6    & 11    &  8    & 14.35 & 11.48 & 10.45 &  8.71 & 11.61 & 12.75 &  9.56 & 14.88 &  9.23 & 13.84 & 13.84 & 0.801 \\
\textnormal{flt2-prod3} & 5     &  6    &  1    & 11    &  2.87 &  3.59 &  6.88 &  3.44 &  5.41 &  6.22 &  4.67 &  7.26 &  4.50 &  6.75 &  6.75 & 0.391 \\
\textnormal{flt2-prod4} & 4     &  4    &  4    &  1    &  4.31 &  2.87 &  1.47 &  2.46 &  4.42 &  3.43 &  2.29 &  3.43 &  2.68 &  4.02 &  4.02 & 0.233 \\
\textnormal{flt2-prod5} & 0     &  2    &  0    &  0    &  0.72 &  0.00 &  0.49 &  1.47 &  0.00 &  1.14 &  1.14 &  1.91 &  0.63 &  0.95 &  0.95 & 0.055 \\ \hline
\textnormal{flt3-prod1} & 20    & 30    & 22    & 28    & 30.07 & 24.25 & 26.08 & 21.73 & 28.98 & 31.82 & 23.87 & 37.12 & 23.02 & 34.54 & 34.54 & 2.000 \\
\textnormal{flt3-prod2} & 22    & 12    & 22    & 16    & 28.71 & 22.97 & 20.90 & 17.42 & 23.22 & 25.50 & 19.13 & 29.75 & 18.45 & 27.68 & 27.68 & 1.603 \\
\textnormal{flt3-prod3} & 10    & 12    &  2    & 22    &  5.74 &  7.18 & 13.76 &  6.88 & 10.81 & 12.44 &  9.33 & 14.52 &  9.00 & 13.51 & 13.51 & 0.782 \\
\textnormal{flt3-prod4} & 8    &  8    &  8    &  2    &  8.61 &  5.74 &  2.95 &  4.92 &  8.85 &  6.86 &  4.57 &  6.86 &  5.36 &  8.04 &  8.04 & 0.465 \\
\textnormal{flt3-prod5} & 0    &  4    &  0    &  0    &  1.44 &  0.00 &  0.98 &  2.95 &  0.00 &  2.29 &  2.29 &  3.81 &  1.26 &  1.89 &  1.89 & 0.110 \\ \hline
$\mathbf{\lambda}_t$    & 128.57& 141.43& 115.71& 145.71& 159.79& 128.86& 138.58& 115.49& 153.98& 169.10& 126.83& 197.29& 122.35& 183.53& 183.53&       \\ \hline
\end{tabular}}
\caption{Schedule change example: estimated demand and parameters ($\alpha=0$)}
\label{table:sch_ex1_sol}
\end{table}

For flights 1 and 2 we get the same estimated primary demand and preference weights as in Table \ref{table:vvrr_sol_alpha0}, since we duplicated that example over the time horizon. For flight 3 we estimate twice the primary demand and preference weights, which was expected, since we artificially doubled the numbers. The method is consistent, and can handle non-homogeneous product set in a mathematically formal way.

It is interesting to note here that the EM algorithm of \citet{vvrr12} cannot distinguish between product $i$ not being available for sale ($i \notin S_t$) as opposed to not being part of the product set ($i \notin I_t$). Because of this, naive application of the EM algorithm would estimate primary demand for non-existing products of flight 1 and 2, grossly overestimating the demand. The total estimated arrival rate is $\sum_{i=1}^{30} \lambda_t = 5324.10$, while using the extended EM algorithm we get $\sum_{i=1}^{30} \lambda_t = 4421.53$.

It is also interesting to mention that using the extended EM algorithm with $\alpha = 1$ we get $\sum_{i=1}^{30} \lambda_t = 2365.71$, and just like before, the total primary demand at time period $t$ is equal to the observed purchases. It is easy to show in general that in case $\alpha=1$ it follows that $\lambda_t = \left( \sum_{i=1}^n z_{it} \right)/s$. If we apply the algorithm with other values of $\alpha \in [0,1]$ we observe a linear decrease of total estimated arrival rate as a function of $\alpha$. The model induced market share of the retailer increases as $\alpha$ increases which decreases the estimated demand.

Algorithm \ref{algo_mfem_ext} is a simple but yet powerful extension of the EM algorithm which can natively handle non-homogenous product set, and can control the availability of the outside alternative by a simple parameter. The price of this extension is that in the M-step we need to solve for the roots of a system of nonlinear equations, instead of having a closed form solution. However, we can use a simple fixed point iteration as the M-step.

\section{Constrained optimization}
\label{sec:constrained_opt}

In this section we develop algorithms to solve the estimation problem with constrained arrival rates (Section \ref{sec:constrained_parameters}). We will extend the minorization-maximization (MM) algorithm of \citet{tarek2016} and also present a solution utilizing the Frank-Wolfe algorithm. We will also incorporate partial availability, non-homogeneous product set, and the ability to control the availability of outside alternative into the model, and show how to estimate the parameters using iterative algorithms.

Consider the incomplete log-likelihood
\footnotesize
\be
\label{eq:loglik_mm}
l_I(\bm{v}, \bm{\lambda}) = \sum_{t=1}^{T} \left[ m_t \log\left( \frac{\lambda_t}{v_{0t} + \sum_{i \in S_t} v_i \cdot o_{it}} \right) - \lambda_t \frac{\sum_{i \in S_t} v_i \cdot o_{it}}{v_{0t} + \sum_{i \in S_t} v_i \cdot o_{it}} +\sum_{i \in S_t} z_{it}\log(v_i \cdot o_{it}) \right]
\ee
\normalsize
which is an extension of the log-likelihood in \citet{tarek2016} with partial availability, as discussed in Section \ref{sec:extend_attr_model}. In case $o_{it} = 1$ when $i \in S_t$, the log-likelihood simplifies to the one considering only fully open and closed assortments. Note that we are using preference weights $\bm v$ in the model, but we could express the model in the utility space by using $\bm v = \exp(\bm \beta)$.

The constrained optimization problem we need to solve is given by
\be
\label{eq:loglik_constr}
\max_{\bm v, \bm \lambda} && l_I(\bm{v}, \bm{\lambda}) \\ \nonumber
\textrm{s.t.} && \\ \nonumber
&& \lambda_t \leq L_t,~t=1,\ldots,T
\ee
where $L_t$ is an upper bound on the arrival rate at time $t$. The motivation behind putting an upper bound on the arrival rates was explained in Section \ref{sec:constrained_parameters}, where we discussed two specific ways of constraining the arrival rates. To solve the constrained maximum likelihood estimation problem, we follow the idea in \citet{tarek2016}. We express the optimization problem as a function of $\bm v$ by applying part of the Karush–Kuhn–Tucker (KKT) conditions on the Lagrangian function, removing $\bm{\lambda}$ from the problem. The Lagrangian function of (\ref{eq:loglik_constr}) becomes
\be
{\mathcal{L}}(\bm{v}, \bm{\lambda},\bm{\mu}) = l_I(\bm{v}, \bm{\lambda}) - \sum_{t=1}^{T}\frac{\mu_t}{v_{0t} + \sum_{i \in S_t} v_i \cdot o_{it}}(\lambda_t-L_t)
\ee
Note that we used a simple re-scaling ($\frac{1}{v_{0t} + \sum_{i \in S_t} v_i \cdot o_{it}}$ for time $t$) of the Lagrange multipliers $\bm \mu$. After we apply the KKT conditions to remove $\bm \lambda$ from the problem, and some algebra (see Appendix), the problem becomes
\be
\label{eq:loglik_constr_reduced}
\max_{\bm v} &&\sum_{i=1}^n K_i\log(v_i) - \sum_{t\not\in\mathcal{B}(\bm{v})} m_t \log\left(\sum_{i \in S_t} v_i \cdot o_{it} \right) - \\ \nonumber
&& \sum_{t\in\mathcal{B}(\bm{v})} m_t \log\left({v_{0t} + \sum_{i \in S_t} v_i \cdot o_{it}}\right) 
- \sum_{t\in\mathcal{B}(\bm{v})} L_t \frac{\sum_{i \in S_t} v_i \cdot o_{it}}{v_{0t} + \sum_{i \in S_t} v_i \cdot o_{it}}+C_3\left(\mathcal{B}(\bm{v})\right)
\ee
where $C_3\left(\mathcal{B}(\bm{v})\right)$ is defined in $(\ref{eq:c_v})$, $K_i = \sum_{t=1}^T z_{it}$ and $\mathcal{B}(\bm{v})$ represents the set of time periods where the upper bound constraints on the arrival rates are binding, that is
\be
\mathcal{B}(\bm{v})= \left\{t \bigg| L_t < m_t\frac{v_{0t} + \sum_{i \in S_t} v_i \cdot o_{it}}{\sum_{i \in S_t} v_i \cdot o_{it}}\right\}
\ee

Through the definitions of $v_{0t}$ and market share constraint we can incorporate non-homoge\-neous product set in the model, and the ability to control the availability of outside alternative. These details were discussed in Sections \ref{sec:nonhomo_prod} and \ref{sec:share_nopurchase}. We will consider two formulations here. In the first one we use market share constraint
\be
\label{eq:constraints_mm}
\sum_{t=1}^T \left[ (1-\alpha) \sum_{i \in I_t} v_i + \alpha \sum_{i \in S_t} v_i \cdot o_{it} \right] = \tilde{s} \sum_{t=1}^T v_{0t}
\ee
which is equivalent to the aggregate constraint (\ref{eq:share_constr_agg}) discussed in Section \ref{sec:share_nopurchase}, with using $\tilde{s} = s/(1-s)$. We assume that OA preference weights $v_{0t}$ are known, and we can use standard scaling $v_{0t}=1$. In Section \ref{sec:extended_mm} we will show how to use the MM algorithm to solve this problem.

In the second formulation we incorporate market share constraint (\ref{eq:share_constr_alpha}) into the objective function and constrain the preference weights by using
\be
\label{eq:constraints_fw}
v_{0t} &=& r \left[(1-\alpha)\sum_{i \in I_t} v_i + \alpha \sum_{i \in S_t} v_i \cdot o_{it} \right] \\ \nonumber
\sum_{\forall i} v_i &=& 1
\ee
This first equation removes $v_{0t}$ from the problem, while the second equation avoids having multiple solutions by rescaling $\bm v$. In Section \ref{sec:frank_wolfe} we will show how the Frank-Wolfe method will lead to a simple coordinate descent algorithm to solve this problem.

\subsection{Solution using MM algorithm}
\label{sec:extended_mm}

In this section we present a solution to the optimization problem (\ref{eq:loglik_constr}) with constraint (\ref{eq:constraints_mm}) using the MM algorithm, building on \citet{tarek2016}. The idea behind MM algorithms is to find a surrogate function that minorizes the original objective function, maximize the surrogate function, and continue this iteratively. For more information on the MM algorithm, in general, see \citet{hunter2000, hunter2004}.

After removing $\bm \lambda$ from the problem using the KKT conditions, we arrive to (\ref{eq:loglik_constr_reduced}). If we rearrange the last term of the objective to aid the minorization, the optimization problem becomes
\be
\label{eq:opt_problem_mm}
\max_{\bm v} &&\sum_{i=1}^n K_i\log(v_i) - \sum_{t\not\in\mathcal{B}(\bm{v})} m_t \log\left(\sum_{i \in S_t} v_i \cdot o_{it} \right) - \\ \nonumber
&& \sum_{t\in\mathcal{B}(\bm{v})} m_t \log\left({v_{0t} + \sum_{i \in S_t} v_i \cdot o_{it}}\right) + \sum_{t\in\mathcal{B}(\bm{v})} L_t \frac{v_{0t}}{v_{0t} + \sum_{i \in S_t} v_i \cdot o_{it}} +C_3\left(\mathcal{B}(\bm{v})\right)\\ \nonumber
\textrm{s.t.}&& \\ \nonumber
&& \sum_{t=1}^T \left[ (1-\alpha) \sum_{i \in I_t} v_i + \alpha \sum_{i \in S_t} v_i \cdot o_{it} \right] = \tilde{s} \sum_{t=1}^T v_{0t}
\ee
where $v_{0t}$ are known. A common technique to find a minorizer is to use supporting hyperplanes \citep{hunter2004}. Since the second, third, and fourth terms in (\ref{eq:opt_problem_mm}) are convex, we can use first-order Taylor approximation to the convex functions $-\log(x)$ and $\frac{1}{x}$, that is, for all $x,y>0$
\ba
-\log(y) &\geq& -\log(x)-\frac{1}{x}(y-x) \\
\frac{1}{y} &\geq& \frac{1}{x} - \frac{1}{x^2}(y-x)
\ea
In our specific case we will use
\ba
-\log(y) &\geq& 1-\log(x) - y/x \\
\frac{1}{y} &\geq& - \frac{y}{x^2}+\frac{2}{x}
\ea
with $y = \sum_{i \in S_t} v_i \cdot o_{it}$ and $x = \sum_{i \in S_t} v^{(k)}_i \cdot o_{it}$ or $y = v_{0t}+\sum_{i \in S_t} v_i \cdot o_{it}$ and $x = v_{0t} + \sum_{i \in S_t} v^{(k)}_i \cdot o_{it}$, where $\bm{v}^{(k)}$ is the value of $\bm{v}$ at iteration $k$. Therefore, the minorizer function becomes
\be
g(\bm{v}|\bm{v}^{(k)}) &=& \sum_{i=1}^n K_i \log(v_i)-\sum_{t\not\in\mathcal{B}(\bm{v}^{(k)})}m_t\frac{\sum_{i \in S_t} v_i \cdot o_{it}}{\sum_{i \in S_t} v_i^{(k)} \cdot o_{it}} \\ \nonumber
&&  -\sum_{t\in\mathcal{B}(\bm{v}^{(k)})}m_t\frac{\sum_{i \in S_t} v_i \cdot o_{it}}{v_{0t} + \sum_{i \in S_t} v_i^{(k)} \cdot o_{it}}
-\sum_{t\in\mathcal{B}(\bm{v}^{(k)})}\frac{L_t v_{0t}(\sum_{i \in S_t} v_i \cdot o_{it})}{\left(v_{0t} + \sum_{i \in S_t} v_i^{(k)} \cdot o_{it}\right)^2} + C_0
\ee
where $C_0$ contains the constant terms independent of $\bm{v}$. Although the minorizer $g(\bm{v}|\bm{v}^{(k)})$ is developed locally for fixed $\bm{v}^{(k)}$, it can be shown to be globally dominated by
the ${\mathcal{L}}(\bm{v}, \bm{\lambda},\bm{\mu})$. Computational results show that it is actually a fairly tight minorizer as well. We now need to solve a single market share
constraint optimization problem:
\ba
\max_{\bm v} && \sum_{i=1}^n K_i \log(v_i)-\sum_{t\not\in\mathcal{B}(\bm{v}^{(k)})}m_t\frac{\sum_{i \in S_t} v_i \cdot o_{it}}{\sum_{i \in S_t} v_i^{(k)} \cdot o_{it}} \\
&&  -\sum_{t\in\mathcal{B}(\bm{v}^{(k)})}m_t\frac{\sum_{i \in S_t} v_i \cdot o_{it}}{v_{0t} + \sum_{i \in S_t} v_i^{(k)} \cdot o_{it}}
-\sum_{t\in\mathcal{B}(\bm{v}^{(k)})}\frac{L_t v_{0t}(\sum_{i \in S_t} v_i \cdot o_{it})}{\left(v_{0t} + \sum_{i \in S_t} v_i^{(k)} \cdot o_{it}\right)^2} \\
\textrm{s.t.}&& \\
&& \sum_{t=1}^T \left[ (1-\alpha) \sum_{i \in I_t} v_i + \alpha \sum_{i \in S_t} v_i \cdot o_{it} \right] = \tilde{s} \sum_{t=1}^T v_{0t}
\ea
Since the objective function is concave, the constraint is convex, and we are maximizing, there exists a single scalar $\eta$ such that the first order optimality condition holds for Lagrangian
\ba
\mathcal{L}(\bm v, \eta) = g(\bm{v}|\bm{v^{(k)}}) - \eta \left( \sum_{t=1}^T \left[ (1-\alpha) \sum_{i \in I_t} v_i + \alpha \sum_{i \in S_t} v_i \cdot o_{it} \right] - \tilde{s} \sum_{t=1}^T v_{0t} \right)
\ea
The first order optimality condition is
\small
\be
\label{eq:Kj}
K_j &=& v_j \left( \sum_{t\not\in\mathcal{B}(\bm{v}^{(k)})} m_t \frac{\mathbbm{1}(j \in S_t)\cdot o_{jt}}{\sum_{i \in S_t} v_i^{(k)}\cdot o_{it}} + \sum_{t\in\mathcal{B}(\bm{v}^{(k)})}m_t \frac{ \mathbbm{1}(j \in S_t)\cdot o_{jt}}{v_{0t} + \sum_{i \in S_t} v_i^{(k)}\cdot o_{it}} + \right. \\ \nonumber
&& \left. \sum_{t\in\mathcal{B}(\bm{v}^{(k)})} L_t \frac{v_{0t} \mathbbm{1}(j \in S_t)\cdot o_{jt}}{\left( v_{0t} + \sum_{i \in S_t} v_i^{(k)}\cdot o_{it}\right)^2} + \eta \left[ (1-\alpha)\sum_{t=1}^T \mathbbm{1}(j \in I_t) + \alpha \sum_{t=1}^T \mathbbm{1}(j \in S_t)\cdot o_{jt} \right] \right)
\ee
\normalsize
which leads to the MM update of $\bm v$ summarized in Algorithm \ref{alg:mm}. In the update step we use Newton's method to find $\eta$ in every MM iteration, summarized in Algorithm \ref{alg:newton}. For more details, please see Appendix.

\begin{algorithm}[!htbp]
\footnotesize
\caption{MM algorithm to estimate $(\bm \lambda, \bm v)$ in (\ref{eq:loglik_constr}) with market share constraint (\ref{eq:constraints_mm})}
\label{alg:mm}
\begin{algorithmic}[1]
\State Input: $\left\{\left(\bm{z}_t,\bm{o}_t,S_t,I_t,v_{0t},L_t \right)\right\}_{t=1}^T$, $s$, $\alpha$
\State Let $K_j = \sum_{t=1}^T z_{jt}$, $n_j = \sum_{t=1}^T \mathbbm{1}(j \in I_t)$ and $o_j = \sum_{t=1}^T o_{jt},~j=1,\ldots,n$
\State Let $m_t = \sum_{j=1}^n z_{jt},~t=1,\ldots,T$ and $\tilde{s} = s/(1-s)$
\State Let $k=0$, initialize $\bm{v}^{(0)}$
\While{Stopping criteria is not satisfied}
    \State $\mathcal{B}(\bm{v}^{(k)})= \left\{t \bigg| L_t < m_t\frac{v_{0t} + \sum_{i \in S_t} v_i^{(k)} \cdot o_{it}}{\sum_{i \in S_t} v_i^{(k)} \cdot o_{it}}\right\}$
    \For{$j = 1,\ldots,n$}
        \State $A_j = \sum_{t \notin \mathcal{B}(\bm{v}^{(k)})} m_t \frac{\mathbbm{1}(j \in S_t)\cdot o_{jt}}{\sum_{i \in S_t} v_i^{(k)}\cdot o_{it}}
        + \sum_{t \in \mathcal{B}(\bm{v}^{(k)})} m_t \frac{\mathbbm{1}(j \in S_t)\cdot o_{jt}}{v_{0t} + \sum_{i \in S_t} v_i^{(k)}\cdot o_{it}} + \sum_{t \in \mathcal{B}(\bm{v}^{(k)})} L_t \frac{v_{0t} \mathbbm{1}(j \in S_t)\cdot o_{jt}}{\left(v_{0t} + \sum_{i \in S_t} v_i^{(k)}\cdot o_{it}\right)^2}$
    \EndFor
    \State Find $\eta$ using \textbf{Algorithm \ref{alg:newton}}
    \For{$j = 1,\ldots,n$}
        \State $v_j^{(k+1)} = \frac{K_j}{A_j + \eta \left[ (1-\alpha) n_j + \alpha o_j \right]}$
    \EndFor
    \State $k = k+1$
\EndWhile
\For{$t = 1,\ldots,T$}
    \State $\lambda_t = \min \left( L_t, m_t\frac{v_{0t} + \sum_{i \in S_t} v_i^{(k)} \cdot o_{it}}{\sum_{i \in S_t} v_i^{(k)} \cdot o_{it}} \right)$
\EndFor
\end{algorithmic}
\end{algorithm}


\begin{algorithm}[!htbp]
\footnotesize
\caption{Newton's method to find $\eta$}
\label{alg:newton}
\begin{algorithmic}[1]
\State Let $k=0$, $\eta^0 = 0$, and $f^0 = 1$
\While{$f^k > \epsilon $}
    \State $f^k = \sum_{j=1}^n \frac{K_j \left[ (1-\alpha) n_j + \alpha o_j \right]}{A_j + \eta^k \left[ (1-\alpha) n_j + \alpha o_j \right]} - \tilde{s} \sum_{t=1}^T v_{0t}$
    \State $g^k = - \sum_{j=1}^n \frac{K_j \left[ (1-\alpha) n_j + \alpha o_j \right]^2}{\left(A_j + \eta^k \left[ (1-\alpha) n_j + \alpha o_j \right]\right)^2}$
    \State $\eta^{k+1} = \eta^{k} - \frac{f^k}{g^k}$
    \State $k = k+1$
\EndWhile
\end{algorithmic}
\end{algorithm}

\subsection{Solution using Frank-Wolfe method}
\label{sec:frank_wolfe}

In this section we present a solution to the optimization problem (\ref{eq:loglik_constr}) with constraint (\ref{eq:constraints_fw}) using the Frank-Wolfe algorithm. The Frank-Wolfe, or conditional gradient method is an iterative optimization algorithm for constrained convex optimization, where in each iteration, it considers a linear approximation of the objective function, and moves towards a minimizer of this linear function. For more information on the Frank-Wolfe algorithm, see \citet{frankwolfe1956} and \citet{bertsekas1999}.

After removing $\bm \lambda$ from the problem using the KKT conditions, we arrive to (\ref{eq:loglik_constr_reduced}). Therefore, we need to solve the constrained optimization problem
\be
\label{eq:opt_problem_fw}
\max_{\bm v} &&\sum_{i=1}^n K_i\log(v_i) - \sum_{t\not\in\mathcal{B}(\bm{v})} m_t \log\left(\sum_{i \in S_t} v_i \cdot o_{it} \right) - \nonumber\\
&& \sum_{t\in\mathcal{B}(\bm{v})} m_t \log\left({v_{0t} + \sum_{i \in S_t} v_i \cdot o_{it}}\right) - \sum_{t\in\mathcal{B}(\bm{v})} L_t \frac{\sum_{i \in S_t} v_i \cdot o_{it}}{v_{0t} + \sum_{i \in S_t} v_i \cdot o_{it}} +C_3\left(\mathcal{B}(\bm{v})\right)\nonumber\\
\textrm{s.t.}&& \nonumber \\
&&\sum_{i=1}^n v_i = 1
\ee
where we plug in $v_{0t} = r \left[(1-\alpha)\sum_{i \in I_t} v_i + \alpha \sum_{i \in S_t} v_i \cdot o_{it} \right]$ to include the market share constraints into the objective function.

The first step of the Frank-Wolfe algorithm is the direction finding subproblem, which becomes
\be
\label{eq:fw_step1}
\max_{\bm y} && \nabla f\left(\bm{v}^{(k)}\right)^T \bm{y}  \\ \nonumber
\textrm{s.t.} && \\ \nonumber
&&\sum_{i=1}^n y_i = 1
\ee
where $\nabla f\left(\bm{v}^{(k)}\right)$ is the gradient vector of the objective function (\ref{eq:opt_problem_fw}) evaluated at the solution of iteration $k$. The elements of $\nabla f\left(\bm{v}\right)$ are calculated as
\small
\ba
\frac{\partial f}{\partial v_j} = && \frac{K_j}{v_j} - \sum_{t\not\in\mathcal{B}(\bm{v})} m_t \frac{\mathbbm{1}(j \in S_t) o_{jt}}{\sum_{i \in S_t} v_i \cdot o_{it}} - \\ && \sum_{t\in\mathcal{B}(\bm{v})} m_t \frac{ \mathbbm{1}(j \in S_t) o_{jt} (r\alpha+1) } {v_{0t} + \sum_{i \in S_t} v_i \cdot o_{it}} - \sum_{t\in\mathcal{B}(\bm{v})} L_t \frac{ \mathbbm{1}(j \in S_t) o_{jt} r(1-\alpha) (\sum_{i \in I_t} v_i)}{(v_{0t} + \sum_{i \in S_t} v_i \cdot o_{it})^2} - \\
&& \sum_{t\in\mathcal{B}(\bm{v})} m_t \frac{ \mathbbm{1}(j \in I_t) r(1-\alpha) }{v_{0t} + \sum_{i \in S_t} v_i \cdot o_{it}} + \sum_{t\in\mathcal{B}(\bm{v})} L_t \frac{ \mathbbm{1}(j \in I_t) r(1-\alpha)(\sum_{i \in S_t} v_i \cdot o_{it}) }{(v_{0t} + \sum_{i \in S_t} v_i \cdot o_{it})^2} \\
\textnormal{where} \\
&& v_{0t} = r \left[(1-\alpha)\sum_{i \in I_t} v_i + \alpha \sum_{i \in S_t} v_i \cdot o_{it} \right]
\ea
\normalsize
For detailed derivation and computational formula, please see Appendix. The direction finding subproblem in (\ref{eq:fw_step1}) is a fractional knapsack problem \citep{korte2012}, which solution becomes
\ba
y_i =
\begin{cases}
    1, & \text{if } i = \arg \max \left\{ \left. \frac{\partial f}{\partial v_i}\right|_{\bm{v} = \bm{v}^{(k)}},~i=1,\ldots,n \right\} \\
    0,              & \text{otherwise}
\end{cases}
\ea
We take a step in the direction of the maximum element of the gradient, only changing that variable in the current step. The Frank-Wolfe algorithm reduces to a coordinate descent algorithm, successively maximizing along coordinate directions determined by the largest value of the gradient vector. The update step of Frank-Wolfe is
\ba
\bm{v}^{(k+1)} = \bm{v}^{(k)} + \gamma_k \left( \bm{y}-\bm{v}^{(k)} \right)
\ea
which simplifies to
\ba
v_i^{(k+1)} =
\begin{cases}
    \gamma_k + (1-\gamma_k) v_i^{(k)}, & \text{if } i = \arg \max \left\{ \left. \frac{\partial f}{\partial v_i}\right|_{\bm{v} = \bm{v}^{(k)}},~i=1,\ldots,n \right\} \\
    v_i^{(k)},              & \text{otherwise}
\end{cases}
\ea
For step size $\gamma_k$ we can use the default choice $\gamma_k = 2/(k+2)$ or perform a line search to find $\gamma_k$ that minimizes $f \left( \bm{v}^{(k)} + \gamma_k \left( \bm{s} - \bm{v}^{(k)} \right) \right)$ subject to $0 \leq \gamma_k \leq 1$. In practice, we implemented a backtracking linesearch using the Armijo's rule \citep{nocedal2006}. The Frank-Wolfe algorithm is summarized in Algorithm \ref{alg:fw}, while the backtracking linesearch is presented in Algorithm \ref{alg:armijo}.

\begin{algorithm}[!htbp]
\footnotesize
\caption{Frank-Wolfe algorithm to estimate $(\bm \lambda, \bm v)$ in (\ref{eq:loglik_constr}) with constraint (\ref{eq:constraints_fw})}
\label{alg:fw}
\begin{algorithmic}[1]
\State Input: $\left\{\left(\bm{z}_t,\bm{o}_t,S_t,I_t,L_t \right)\right\}_{t=1}^T$, $s$, $\alpha$
\State Let $K_j = \sum_{t=1}^T z_{jt},~j=1,\ldots,n$
\State Let $m_t = \sum_{j=1}^n z_{jt},~t=1,\ldots,T$ and $r = (1-s)/s$
\State Let $k=0$, initialize $\bm{v}^{(0)}$
\While{Stopping criteria is not satisfied}
    \For{$t = 1,\ldots,T$}
        \State $v_{0t}^{(k)} = r \left[(1-\alpha)\sum_{i \in I_t} v_i^{(k)} + \alpha \sum_{i \in S_t} v_i^{(k)} \cdot o_{it} \right]$
    \EndFor
    \State $\mathcal{B}(\bm{v}^{(k)})= \left\{t \bigg| L_t < m_t\frac{v_{0t}^{(k)} + \sum_{i \in S_t} v_i^{(k)} \cdot o_{it}}{\sum_{i \in S_t} v_i^{(k)} \cdot o_{it}}\right\}$
    \For{$j = 1,\ldots,n$}
        \State
            \ba
            g_j = && \frac{K_j}{v_j^{(k)}} - \sum_{t\not\in\mathcal{B}(\bm{v}^{(k)})} m_t \frac{\mathbbm{1}(j \in S_t) o_{jt}}{\sum_{i \in S_t} v_i^{(k)} \cdot o_{it}} - \\ && \sum_{t\in\mathcal{B}(\bm{v}^{(k)})} m_t \frac{ \mathbbm{1}(j \in S_t) o_{jt} (r\alpha+1) } {v_{0t}^{(k)} + \sum_{i \in S_t} v_i^{(k)} \cdot o_{it}} - \sum_{t\in\mathcal{B}(\bm{v}^{(k)})} L_t \frac{ \mathbbm{1}(j \in S_t) o_{jt} r(1-\alpha) \left(\sum_{i \in I_t} v_i^{(k)}\right)}{\left(v_{0t}^{(k)} + \sum_{i \in S_t} v_i^{(k)} \cdot o_{it}\right)^2} - \\
            && \sum_{t\in\mathcal{B}(\bm{v}^{(k)})} m_t \frac{ \mathbbm{1}(j \in I_t) r(1-\alpha) }{v_{0t}^{(k)} + \sum_{i \in S_t} v_i^{(k)} \cdot o_{it}} + \sum_{t\in\mathcal{B}(\bm{v}^{(k)})} L_t \frac{ \mathbbm{1}(j \in I_t) r(1-\alpha) \left(\sum_{i \in S_t} v_i^{(k)} \cdot o_{it}\right) }{\left(v_{0t}^{(k)} + \sum_{i \in S_t} v_i^{(k)} \cdot o_{it}\right)^2}
            \ea
    \EndFor
    \State $l = \arg \max \left\{ g_j,~j=1,\ldots,n \right\}$ \Comment Find direction
    \State $\bm{y} = \bm{e}_l$ \Comment $\bm{e}_i$ is $i$th unit vector
    \State Use $\gamma_k = 2/(k+2)$ or find $\gamma_k$ using Algorithm \ref{alg:armijo} \Comment Find step size
    \State $\bm{v}^{(k+1)} = \bm{v}^{(k)} + \gamma_k \left( \bm{y}-\bm{v}^{(k)} \right)$ \Comment Frank-Wolfe update
    \State $k = k+1$
\EndWhile
\For{$t = 1,\ldots,T$}
    \State $v_{0t}^{(k)} = r \left[(1-\alpha)\sum_{i \in I_t} v_i^{(k)} + \alpha \sum_{i \in S_t} v_i^{(k)} \cdot o_{it} \right]$
    \State $\lambda_t = \min \left( L_t, m_t\frac{v_{0t}^{(k)} + \sum_{i \in S_t} v_i^{(k)} \cdot o_{it}}{\sum_{i \in S_t} v_i^{(k)} \cdot o_{it}} \right)$
\EndFor
\end{algorithmic}
\end{algorithm}


\begin{algorithm}[!htbp]
\footnotesize
\caption{Backtracking linesearch by Armijo's rule to find step size}
\label{alg:armijo}
\begin{algorithmic}[1]
\State Initialize $\beta$, $\tau$, $\gamma^{(0)}$ \Comment e.g. $\beta = 0.001$, $\tau = 0.5$, $\gamma^{(0)} = 1$
\State Let $h=0$, $a=1$, $b=2$
\State Let $\bm s = \bm y - \bm v^{(h)}$
\State Let
        \ba
        f_v = && \sum_{i=1}^n K_i\log(v^{(k)}_i) - \sum_{t\not\in\mathcal{B}(\bm{v}^{(k)})} m_t \log\left(\sum_{i \in S_t} v_i^{(k)} \cdot o_{it} \right) - \\
        && \sum_{t\in\mathcal{B}(\bm{v}^{(k)})} m_t \log\left({v_{0t}^{(k)} + \sum_{i \in S_t} v_i^{(k)} \cdot o_{it}}\right) - \sum_{t\in\mathcal{B}(\bm{v}^{(k)})} L_t \frac{\sum_{i \in S_t} v_i^{(k)} \cdot o_{it}}{v_{0t}^{(k)} + \sum_{i \in S_t} v_i^{(k)} \cdot o_{it}}
        +C_3\left(\mathcal{B}(\bm{v}^{(k)})\right)
        \ea
\State Let $g_vs = \sum_{i=1}^n g_i s_i$
\While{$a < b$}
    \State $\bm{w} = \bm{v}^{(k)} + \gamma^{(h+1)}\bm{s}$
    \State $w_{0t} = r \left[(1-\alpha)\sum_{i \in I_t} w_i + \alpha \sum_{i \in S_t} w_i \cdot o_{it} \right],~t=1,\ldots,T$
    \State
        \ba
         a &=& \sum_{i=1}^n K_i\log(w_i) - \sum_{t\not\in\mathcal{B}(\bm{w})} m_t \log\left(\sum_{i \in S_t} w_i \cdot o_{it} \right) - \\
         && \sum_{t\in\mathcal{B}(\bm{w})} m_t \log\left({w_{0t} + \sum_{i \in S_t} w_i \cdot o_{it}}\right) - \sum_{t\in\mathcal{B}(\bm{w})} L_t \frac{\sum_{i \in S_t} w_i \cdot o_{it}}{w_{0t} + \sum_{i \in S_t} w_i \cdot o_{it}}+C_3\left(\mathcal{B}(\bm{w})\right)
         \ea
    \State $b = f_v + \gamma^{(h+1)} \beta g_vs$
    \State $\gamma^{(h+1)} = \tau \gamma^{(h)}$
    \State $h = h+1$
\EndWhile
\end{algorithmic}
\end{algorithm}

\subsection{Example: non-homogeneous product set with constraint}
To demonstrate Algorithms \ref{alg:mm} and \ref{alg:fw}, let us revisit the example presented in Table \ref{table:sch_ex1} by adding constraints to the problem. We constrain the arrival rates to be less than twice the observed sales, that is $L_t = 2 m_t$. We use $\alpha=0$, assuming that the OA product is always available, and set $v_{0t} = 1$ for the MM algorithm. Note that $o_t = S_t$, since we have fully open and closed assortments. The results, using the MM algorithm, are presented in Table \ref{table:sch_ex1_mmsol}.

\begin{table}[!htbp]
\centering
\resizebox{0.95\linewidth}{!}{
\begin{tabular}[!htbp]{|c||c|c|c|c|c|c|c|c|c|c|c|c|c|c|c||c|}
	\hline
\multirow{2}{*}{\bf Estimates} & \multicolumn{15}{c||}{\bf Period} & \multirow{2}{*}{$\mathbf{v}_i$} \\ \cline{2-16}
& \bf 1 & \bf 2 & \bf 3 & \bf 4 & \bf 5 & \bf 6 & \bf 7 & \bf 8 & \bf 9 & \bf 10 & \bf 11 & \bf 12 & \bf 13 & \bf 14 & \bf 15 & \\ \hline
\textnormal{flt1-prod1} & 10.41 & 11.45 &  9.37 & 11.80 & 12.47 & 10.06 &  8.74 &  7.29 &  9.72 &  5.83 &  4.37 &  6.80 &  0.97 &  1.46 &  1.46
                        & 1.000 \\
\textnormal{flt1-prod2} &  9.40 & 10.34 &  8.46 & 10.65 & 11.26 &  9.08 &  7.89 &  6.58 &  8.77 &  5.26 &  3.95 &  6.14 &  0.88 &  1.32 &  1.32
                        & 0.903 \\
\textnormal{flt1-prod3} &  5.11 &  5.62 &  4.60 &  5.79 &  6.12 &  4.94 &  4.29 &  3.58 &  4.77 &  2.86 &  2.15 &  3.34 &  0.48 &  0.72 &  0.72
                        & 0.491 \\
\textnormal{flt1-prod4} &  3.71 &  4.08 &  3.33 &  4.20 &  4.44 &  3.58 &  3.11 &  2.59 &  3.46 &  2.07 &  1.56 &  2.42 &  0.35 &  0.52 &  0.52
                        & 0.356 \\
\textnormal{flt1-prod5} &  1.38 &  1.52 &  1.24 &  1.56 &  1.65 &  1.33 &  1.16 &  0.97 &  1.29 &  0.77 &  0.58 &  0.90 &  0.13 &  0.19 &  0.19
                        & 0.133 \\ \hline
\textnormal{flt2-prod1} & \multicolumn{15}{c|}{\gcl}
                        & 1.000 \\
\textnormal{flt2-prod2} & \multicolumn{15}{c|}{\gcl}
                        & 0.903 \\
\textnormal{flt2-prod3} & \multicolumn{15}{c|}{\gcl}
                        & 0.491 \\
\textnormal{flt2-prod4} & \multicolumn{15}{c|}{\gcl}
                        & 0.356 \\
\textnormal{flt2-prod5} & \multicolumn{15}{c|}{\gcl}
                        & 0.133 \\ \hline
\textnormal{flt3-prod1} & 20.82 & 22.90 & 18.74 & 23.59 & 24.94 & 20.11 & 17.49 & 14.57 & 19.43 & 11.66 &  8.74 & 13.60 &  1.94 &  2.91 &  2.91
                        & 2.000 \\
\textnormal{flt3-prod2} & 18.79 & 20.67 & 16.91 & 21.30 & 22.52 & 18.16 & 15.79 & 13.16 & 17.54 & 10.52 &  7.89 & 12.28 &  1.75 &  2.63 &  2.63
                        & 1.806 \\
\textnormal{flt3-prod3} & 10.22 & 11.24 &  9.20 & 11.58 & 12.24 &  9.87 &  8.58 &  7.15 &  9.54 &  5.72 &  4.29 &  6.68 &  0.95 &  1.43 &  1.43
                        & 0.982 \\
\textnormal{flt3-prod4} &  7.41 &  8.15 &  6.67 &  8.40 &  8.88 &  7.16 &  6.22 &  5.19 &  6.92 &  4.15 &  3.11 &  4.84 &  0.69 &  1.04 &  1.04
                        & 0.712 \\
\textnormal{flt3-prod5} &  2.76 &  3.03 &  2.48 &  3.13 &  3.31 &  2.67 &  2.32 &  1.93 &  2.57 &  1.54 &  1.16 &  1.80 &  0.26 &  0.39 &  0.39
                        & 0.265 \\ \hline
$\mathbf{\lambda}_t$    & 128.57& 141.43& 115.71& 145.71& 154.04& 124.22& 108.00&  90.00& 120.00&  72.00&  54.00&  84.00&  12.00&  18.00&  18.00
                        & \\ \hline
\end{tabular}}
\resizebox{0.95\linewidth}{!}{
\begin{tabular}[!htbp]{|c||c|c|c|c|c|c|c|c|c|c|c|c|c|c|c||c|}
	\hline
\bf \phantom{Estimates} & \bf 16 & \bf 17 & \bf 18 & \bf 19 & \bf 20 & \bf 21 & \bf 22 & \bf 23 & \bf 24 & \bf 25 & \bf 26 & \bf 27 & \bf 28 & \bf 29 & \bf 30 & \phantom{$\mathbf{v}_i$} \\ \hline
\textnormal{flt1-prod1} & \multicolumn{15}{c|}{\gcl}    & 1.000 \\
\textnormal{flt1-prod2} & \multicolumn{15}{c|}{\gcl}    & 0.903 \\
\textnormal{flt1-prod3} & \multicolumn{15}{c|}{\gcl}    & 0.491 \\
\textnormal{flt1-prod4} & \multicolumn{15}{c|}{\gcl}    & 0.356 \\
\textnormal{flt1-prod5} & \multicolumn{15}{c|}{\gcl}    & 0.133 \\ \hline
\textnormal{flt2-prod1} & 10.41 & 11.45 &  9.37 & 11.80 & 12.47 & 10.06 &  8.74 &  7.29 &  9.72 &  5.83 &  4.37 &  6.80 &  0.97 &  1.46 &  1.46 & 1.000 \\
\textnormal{flt2-prod2} &  9.40 & 10.34 &  8.46 & 10.65 & 11.26 &  9.08 &  7.89 &  6.58 &  8.77 &  5.26 &  3.95 &  6.14 &  0.88 &  1.32 &  1.32 & 0.903 \\
\textnormal{flt2-prod3} &  5.11 &  5.62 &  4.60 &  5.79 &  6.12 &  4.94 &  4.29 &  3.58 &  4.77 &  2.86 &  2.15 &  3.34 &  0.48 &  0.72 &  0.72 & 0.491 \\
\textnormal{flt2-prod4} &  3.71 &  4.08 &  3.33 &  4.20 &  4.44 &  3.58 &  3.11 &  2.59 &  3.46 &  2.07 &  1.56 &  2.42 &  0.35 &  0.52 &  0.52 & 0.356 \\
\textnormal{flt2-prod5} &  1.38 &  1.52 &  1.24 &  1.56 &  1.65 &  1.33 &  1.16 &  0.97 &  1.29 &  0.77 &  0.58 &  0.90 &  0.13 &  0.19 &  0.19 & 0.133 \\ \hline
\textnormal{flt3-prod1} & 20.82 & 22.90 & 18.74 & 23.59 & 24.94 & 20.11 & 17.49 & 14.57 & 19.43 & 11.66 &  8.74 & 13.60 &  1.94 &  2.91 &  2.91 & 2.000 \\
\textnormal{flt3-prod2} & 18.79 & 20.67 & 16.91 & 21.30 & 22.52 & 18.16 & 15.79 & 13.16 & 17.54 & 10.52 &  7.89 & 12.28 &  1.75 &  2.63 &  2.63 & 1.806 \\
\textnormal{flt3-prod3} & 10.22 & 11.24 &  9.20 & 11.58 & 12.24 &  9.87 &  8.58 &  7.15 &  9.54 &  5.72 &  4.29 &  6.68 &  0.95 &  1.43 &  1.43 & 0.982 \\
\textnormal{flt3-prod4} &  7.41 &  8.15 &  6.67 &  8.40 &  8.88 &  7.16 &  6.22 &  5.19 &  6.92 &  4.15 &  3.11 &  4.84 &  0.69 &  1.04 &  1.04 & 0.712 \\
\textnormal{flt3-prod5} &  2.76 &  3.03 &  2.48 &  3.13 &  3.31 &  2.67 &  2.32 &  1.93 &  2.57 &  1.54 &  1.16 &  1.80 &  0.26 &  0.39 &  0.39 & 0.265 \\ \hline
$\mathbf{\lambda}_t$    & 128.57& 141.43& 115.71& 145.71& 154.04& 124.22& 108.00&  90.00& 120.00&  72.00&  54.00&  84.00&  12.00&  18.00&  18.00&       \\ \hline
\end{tabular}}
\caption{Schedule change example with constraint: estimated demand and parameters using MM algorithm}
\label{table:sch_ex1_mmsol}
\end{table}

We observe the expected symmetry in the results, that is the estimated demands and preference weights are same for flights 1 and 2, and twice for flight 3. We can see that the estimated value of $\lambda_t$ is equal to the upper bound $L_t$ for time periods 7-15 and 22-30. Using the Frank-Wolfe algorithm, we converge to the same solution. Note that Frank-Wolfe is solving a problem with different constraints, but for this symmetric example and in the case of $\alpha=0$ the two different formulations are equivalent. It is also interesting to note that using a built-in nonlinear optimization routine in R (sequential quadratic programming), we can solve the problem in 141 iterations, taking 2.53 seconds, converging to the same solution. The Frank-Wolfe algorithm with Armijo's rule converges in 267 iterations taking 0.12 second, and the MM algorithm converges in 11 iterations taking only 0.006 seconds. The convergence of the Frank–Wolfe algorithm is sublinear, in general, and using the default step size $\gamma_k = 2/(k+2)$ results in even much slower convergence. Note, however, that the Frank-Wolfe algorithm solves an optimization problem with market share constraint at each time period $t$ substituted into the objective function. For this problem we were not able to develop the MM algorithm due to the difficulty of finding a suitable minorizer. The formulation of the MM algorithm uses an aggregate market share constraint, where we assume that $v_{0t}$ are known. The solutions are not equivalent, in general.

\section{Conclusions and Future Research}
In this paper we discussed some of the practical limitations of the standard choice-based demand models used in the literature to estimate demand from sales transaction data. We presented modifications and extensions of the models and discussed data preprocessing and solution techniques which can be useful for practitioners dealing with sales transaction data. We hope that these discussions could facilitate further methodological progress and even more rigorous theoretical research in this domain. We presented an algorithm to split sales transaction data observed under partial availability, and we extended an EM algorithm for the case we observe a non-homogeneous product set. We developed two iterative optimization algorithms which incorporated partial availability information, non-homogeneous product set, ability to control the availability of outside alternative, and an upper bound on the arrival rates. In one formulation we used market share constraint at each time period, and incorporated them into the objective function through the preference weights of the outside alternative. This formulation was solved using the Frank-Wolfe algorithm, leading to a simple coordinate descent algorithm. We discussed another formulation, which used a single, aggregate market share constraint over the time horizon, and assumed the knowledge of preference weights of the outside alternative. Using this formulation we could develop a very fast, iterative minorization-maximization algorithm building on the work in \citet{tarek2016}.

Future extension of these methods are possible. For instance, after using the sales splitting algorithm, it would be interesting to group the arrival rates in the EM algorithm and avoid having too many parameters to be estimated. Similarly, it would be of practical interest to extend the EM algorithm to the constrained optimization case or introduce other regularization on the parameters. The MM and the Frank-Wolfe algorithms developed in this paper could be extended to include covariates and additional preference weights for the product with lowest fare. The Frank-Wolfe algorithm could be further improved by taking the gradient descent direction instead of coordinate descent. Finally, in practical applications, we often encounter sparse demand distribution with excess number of zeros and heavy tails. A natural extension of the currently used models would be to use a zero-inflated Poisson or negative binomial distribution to describe the customer arrival process.

\newpage
\section{Appendix}
\underline{Using KKT conditions to remove $\bm{\lambda}$} \\ \\
The KKT condition of Lagrangian (\ref{eq:loglik_constr}) are
\ba
\frac{\partial {\mathcal{L}}(\bm{v}, \bm{\lambda},\bm{\mu})}{\partial \lambda_t} = \frac{m_t}{\lambda_t} - \frac{\mu_t + \sum_{i \in S_t} v_i \cdot o_{it}}{v_{0t} + \sum_{i \in S_t} v_i \cdot o_{it}} = 0, ~t=1,\ldots,T \\
\mu_t (\lambda_t-L_t) = 0, ~t=1,\ldots,T \\
\mu_t \geq 0, ~t=1,\ldots,T \\
\lambda_t \leq L_t, ~t=1,\ldots,T
\ea
Expressing $\lambda_t$ and $\mu_t$ in first equation results in
\ba
\lambda_t &=& m_t \frac{v_{0t} + \sum_{i \in S_t} v_i \cdot o_{it}}{\mu_t + \sum_{i \in S_t} v_i \cdot o_{it}} \\
\mu_t &=& \frac{m_t}{\lambda_t}(v_{0t} + \sum_{i \in S_t} v_i \cdot o_{it}) - {\sum_{i \in S_t} v_i \cdot o_{it}}
\ea
Plugging back $\mu_t$, the complementary slackness condition becomes
\ba
\left( \frac{m_t}{\lambda_t}(v_{0t} + \sum_{i \in S_t} v_i \cdot o_{it}) - {\sum_{i \in S_t} v_i \cdot o_{it}} \right) (\lambda_t-L_t) = 0
\ea
Therefore, we have that one of the following conditions should hold:
\ba
\lambda_t^1 &=& L_t \\
\lambda_t^2 &=& m_t\frac{v_{0t} + \sum_{i \in S_t} v_i \cdot o_{it}}{\sum_{i \in S_t} v_i \cdot o_{it}} \\
\ea
So, the partial KKT condition stated earlier reduces to the following conditions: \\
\noindent If $L_t < m_t\frac{v_{0t} + \sum_{i \in S_t} v_i \cdot o_{it}}{\sum_{i \in S_t} v_i \cdot o_{it}}$ then
\ba
\lambda_t &=& L_t \\
\mu_t &=& \frac{m_t}{L_t}(v_{0t} + \sum_{i \in S_t} v_i \cdot o_{it}) - {\sum_{i \in S_t} v_i \cdot o_{it}}\\
\ea
else
\ba
\lambda_t &=& m_t\frac{v_{0t} + \sum_{i \in S_t} v_i \cdot o_{it}}{\sum_{i \in S_t} v_i \cdot o_{it}} \\
\mu_t &=& 0 \\
\ea
Let us define $\mathcal{B}(\bm{v})= \left\{t \bigg| L_t < m_t\frac{v_{0t} + \sum_{i \in S_t} v_i \cdot o_{it}}{\sum_{i \in S_t} v_i \cdot o_{it}}\right\}
$ The Lagrangian function can be simplified to
\ba
{\mathcal{L}}(\bm{v}, \bm{\lambda},\bm{\mu})
&=& l_I(\bm{v}, \bm{\lambda}) - \sum_{t=1}^{T}\left (\frac{\mu_t}{v_{0t} + \sum_{i \in S_t} v_i \cdot o_{it}}\right )(\lambda_t-L_t) \\
&=& \sum_{t=1}^{T} \sum_{i \in S_t} z_{it} \log(v_i \cdot o_{it}) - \sum_{t\not\in\mathcal{B}(\bm{v})} m_t \log\left(\sum_{i \in S_t} v_i \cdot o_{it} \right) + \\
&& \sum_{t\in\mathcal{B}(\bm{v})} m_t \log\left(\frac{L_t}{v_{0t}+\sum_{i \in S_t} v_i \cdot o_{it}}\right) - \sum_{t\in\mathcal{B}(\bm{v})} L_t \frac{\sum_{i \in S_t} v_i \cdot o_{it}}{v_{0t}+\sum_{i \in S_t} v_i \cdot o_{it}} + C_1\left(\mathcal{B}(\bm{v})\right)\\
&=& \sum_{t=1}^{T} \sum_{i \in S_t} z_{it} \log(v_i \cdot o_{it}) - \sum_{t\not\in\mathcal{B}(\bm{v})} m_t \log\left(\sum_{i \in S_t} v_i \cdot o_{it} \right) - \\
&& \sum_{t\in\mathcal{B}(\bm{v})} m_t \log\left(v_{0t}+\sum_{i \in S_t} v_i \cdot o_{it}\right) - \sum_{t\in\mathcal{B}(\bm{v})} L_t \frac{\sum_{i \in S_t} v_i \cdot o_{it}}{v_{0t}+\sum_{i \in S_t} v_i \cdot o_{it}} +C_2\left(\mathcal{B}(\bm{v})\right)\\
&=& \sum_{i=1}^n K_i\log(v_i) - \sum_{t\not\in\mathcal{B}(\bm{v})} m_t \log\left(\sum_{i \in S_t} v_i \cdot o_{it} \right) - \\
&& \sum_{t\in\mathcal{B}(\bm{v})} m_t \log\left({v_{0t} + \sum_{i \in S_t} v_i \cdot o_{it}}\right) - \sum_{t\in\mathcal{B}(\bm{v})} L_t \frac{\sum_{i \in S_t} v_i \cdot o_{it}}{v_{0t}+\sum_{i \in S_t} v_i \cdot o_{it}} \\
&=& \sum_{i=1}^n K_i\log(v_i) - \sum_{t\not\in\mathcal{B}(\bm{v})} m_t \log\left(\sum_{i \in S_t} v_i \cdot o_{it} \right) - \\
&& \sum_{t\in\mathcal{B}(\bm{v})} m_t \log\left({v_{0t} + \sum_{i \in S_t} v_i \cdot o_{it}}\right) + \sum_{t\in\mathcal{B}(\bm{v})} L_t \frac{v_{0t}}{v_{0t} + \sum_{i \in S_t} v_i \cdot o_{it}}+C_3\left(\mathcal{B}(\bm{v})\right)
\ea
where
\be
C_1\left(\mathcal{B}(\bm{v})\right) &=& \sum_{t\not\in\mathcal{B}(\pmb{v})}\left(m_t \log m_t - m_t\right)\nonumber\\
C_2\left(\mathcal{B}(\bm{v})\right) &=& \sum_{t\not\in\mathcal{B}(\pmb{v})}\left(m_t \log m_t - m_t\right)+\sum_{t\in\mathcal{B}(\pmb{v})}m_t \log L_t\nonumber\\
C_3\left(\mathcal{B}(\bm{v})\right) &=& \sum_{t\not\in\mathcal{B}(\pmb{v})}\left(m_t \log m_t - m_t\right)+\sum_{t\in\mathcal{B}(\pmb{v})}\left(m_t \log L_t - L_t\right).\label{eq:c_v}
\ee
Note that
$$
\sum_{t=1}^{T} \sum_{i \in S_t} z_{it} \log(v_i \cdot o_{it}) = \sum_{t=1}^{T} \sum_{i \in S_t} z_{it} \log(v_i) + C_4 = \sum_{i=1}^n K_i\log(v_i) + C_4
$$ for some constant $C_4$, because $z_{it} = 0$ when $i \notin S_t$.

\noindent \underline{Newton's method to find $\eta$} \\ \\
To simplify notation, let us define
\ba
A_j &=& \sum_{t\not\in\mathcal{B}_j(\bm{v}^{(k)})} m_t \frac{\mathbbm{1}(j \in S_t)\cdot o_{jt}}{\sum_{i \in S_t} v_i^{(k)}} + \sum_{t\in\mathcal{B}_j(\bm{v}^{(k)})}m_t \frac{ \mathbbm{1}(j \in S_t)\cdot o_{jt}}{v_{0t} + \sum_{i \in S_t} v_i^{(k)}} + \\
&& \sum_{t\in\mathcal{B}_j(\bm{v}^{(k)})} L_t \frac{v_{0t} \mathbbm{1}(j \in S_t)\cdot o_{jt}}{{(v_{0t} + \sum_{i \in S_t} v_i^{(k)})}^2} \\
n_j &=& \sum_{t=1}^T \mathbbm{1}(j \in I_t) \\
o_j &=& \sum_{t=1}^T \mathbbm{1}(j \in S_t)\cdot o_{jt} = \sum_{t=1}^T o_{jt}
\ea
where $\mathcal{B}_j(\bm{v}^{(k)}) = \mathcal{B}(\bm{v}^{(k)})\cap\{t|j\in S_t\}$, $n_j$ represents the number of times product $j$ is in the offer set $I_t$ over time horizon $T$, and $o_j$ represents the proportion of time product $j$ is available $S_t$ over time horizon $T$. Using the notation above, equation (\ref{eq:Kj}) simplifies to
\ba
K_j &=& v_j A_j + v_j \eta \left[ (1-\alpha) n_j + \alpha o_j \right]
\ea
Expressing $v_j$ in the above equation and plugging it into the market share constraint (\ref{eq:constraints_mm}) leads to
\small
\ba
&& (1-\alpha) \sum_{t=1}^T \sum_{j \in I_t} \frac{K_j}{A_j + \eta \left[ (1-\alpha) n_j + \alpha o_j \right]} + \alpha \sum_{t=1}^T \sum_{j \in S_t} \frac{K_j \cdot o_{jt}}{A_j + \eta \left[ (1-\alpha) n_j + \alpha o_j \right]} = \tilde{s} \sum_{t=1}^T v_{0t}
\ea
\normalsize
and finally to
\ba
\sum_{j=1}^n \frac{K_j \left[ (1-\alpha) n_j + \alpha o_j \right]}{A_j + \eta \left[ (1-\alpha) n_j + \alpha o_j \right]} = \tilde{s} \sum_{t=1}^T v_{0t}
\ea
where we used general equalities
\ba
\sum_{t=1}^T \sum_{j \in I_t} K_j &=& \sum_{j=1}^n n_j K_j \\
\sum_{t=1}^T \sum_{j \in S_t} K_j \cdot o_{jt} &=& \sum_{j=1}^n o_j K_j
\ea
This leads to
\ba
f = \sum_{j=1}^n \frac{K_j \left[ (1-\alpha) n_j + \alpha o_j \right]}{A_j + \eta \left[ (1-\alpha) n_j + \alpha o_j \right]} - \tilde{s} \sum_{t=1}^T v_{0t}
\ea
and
\ba
g = \frac{d f}{d \eta} = - \sum_{j=1}^n \frac{K_j \left[ (1-\alpha) n_j + \alpha o_j \right]^2}{\left(A_j + \eta \left[ (1-\alpha) n_j + \alpha o_j \right]\right)^2}
\ea
and to Newton's method summarized in Algorithm (\ref{alg:newton}).

\noindent \underline{Computation of gradient in Frank-Wolfe algorithm} \\ \\
Here we will derive the gradient vector of the objective function (\ref{eq:opt_problem_fw}), and show a computational formula. The elements of $\nabla f\left(\bm{v}\right)$ can be derived as
\scriptsize
\ba
\frac{\partial f}{\partial v_i} = && \frac{K_i}{v_i} - \sum_{t\not\in\mathcal{B}(\bm{v})} m_t \frac{\mathbbm{1}(i \in S_t) o_{it}}{\sum_{i \in S_t} v_i \cdot o_{it}} -  \sum_{t\in\mathcal{B}(\bm{v})} m_t \frac{ r(1-\alpha)\mathbbm{1}(i \in I_t) + (r\alpha+1) \mathbbm{1}(i \in S_t) o_{it} }{v_{0t} + \sum_{i \in S_t} v_i \cdot o_{it}} - \\
&& \sum_{t\in\mathcal{B}(\bm{v})} L_t \frac{ \mathbbm{1}(i \in S_t) o_{it} (v_{0t} + \sum_{i \in S_t} v_i \cdot o_{it}) - (\sum_{i \in S_t} v_i \cdot o_{it}) (r(1-\alpha)\mathbbm{1}(i \in I_t) + (r\alpha+1) \mathbbm{1}(i \in S_t) o_{it} ) }{(v_{0t} + \sum_{i \in S_t} v_i \cdot o_{it})^2} \\
= && \frac{K_i}{v_i} - \sum_{t\not\in\mathcal{B}(\bm{v})} m_t \frac{\mathbbm{1}(i \in S_t) o_{it}}{\sum_{i \in S_t} v_i \cdot o_{it}} - \sum_{t\in\mathcal{B}(\bm{v})} m_t \frac{ r(1-\alpha)\mathbbm{1}(i \in I_t) + (r\alpha+1) \mathbbm{1}(i \in S_t) o_{it} }{v_{0t} + \sum_{i \in S_t} v_i \cdot o_{it}} - \\
&& \sum_{t\in\mathcal{B}(\bm{v})} L_t \frac{ r(1-\alpha) (\sum_{i \in I_t} v_i) \mathbbm{1}(i \in S_t) o_{it} - r(1-\alpha)(\sum_{i \in S_t} v_i \cdot o_{it})\mathbbm{1}(i \in I_t) }{(v_{0t} + \sum_{i \in S_t} v_i \cdot o_{it})^2} \\
= && \frac{K_i}{v_i} - \sum_{t\not\in\mathcal{B}(\bm{v})} m_t \frac{\mathbbm{1}(i \in S_t) o_{it}}{\sum_{i \in S_t} v_i \cdot o_{it}} - \\ && \sum_{t\in\mathcal{B}(\bm{v})} m_t \frac{ \mathbbm{1}(i \in S_t) o_{it} (r\alpha+1) } {v_{0t} + \sum_{i \in S_t} v_i \cdot o_{it}} - \sum_{t\in\mathcal{B}(\bm{v})} L_t \frac{ \mathbbm{1}(i \in S_t) o_{it} r(1-\alpha) (\sum_{i \in I_t} v_i)}{(v_{0t} + \sum_{i \in S_t} v_i \cdot o_{it})^2} - \\
&& \sum_{t\in\mathcal{B}(\bm{v})} m_t \frac{ \mathbbm{1}(i \in I_t) r(1-\alpha) }{v_{0t} + \sum_{i \in S_t} v_i \cdot o_{it}} + \sum_{t\in\mathcal{B}(\bm{v})} L_t \frac{ \mathbbm{1}(i \in I_t) r(1-\alpha)(\sum_{i \in S_t} v_i \cdot o_{it}) }{(v_{0t} + \sum_{i \in S_t} v_i \cdot o_{it})^2}
\ea
\normalsize
To implement the gradient with matrix-vector operations let us define
\footnotesize
\ba
\bm{K} \in \mathbb{R}^n_{+} &-& K_i \textnormal{ is the total purchases of product $i$ over the selling horizon} \\
\bm{v} \in \mathbb{R}^n_{+} &-& v_i \textnormal{ is the preference weight for product $i$} \\
\bm{v}_{0} \in \mathbb{R}^T_{+} &-& v_{0t} \textnormal{ is the preference weight for outside alternative at time $t$} \\
r \in \mathbb{R}_{+} &-& r = (1-s)/s \textnormal{ , where $s$ is the host market share} \\
\bm{I} \in \{0,1\}^{nxT}  &-& I_{it} = 1 \textnormal{ if product $i$ is in the offer set at time $t$} \\
\bm{S} \in \{0,1\}^{nxT} &-& S_{it} = 1 \textnormal{ if product $i$ is available for sale at time $t$} \\
\bm{O} \in [0,1]^{nxT} &-& o_{it} \textnormal{ is percentage of time product $i$ is available for sale during time period $t$} \\
\bm{m} \in \mathbb{R}^T_{+} &-& m_t \textnormal{ is the total purchases of all products at time $t$} \\
\bm{L} \in \mathbb{R}^T_{+}  &-& L_{t} \textnormal{ is the upper bound on $\lambda_t$} \\
\bm{B} \in \{0,1\}^{T}  &-& B_{t} = 1 \textnormal{ if bound $L_t$ is violated at time $t$} \left(L_t < m_t\frac{v_{0t} + \sum_{i \in S_t} v_i \cdot o_{it}}{\sum_{i \in S_t} v_i \cdot o_{it}}\right) \\
\alpha \in [0,1] &-& \textnormal{parameter to control availability of outside alternative}
\ea
\normalsize
Then
\ba
\nabla f\left(\bm{v}\right) = \bm{K} \oslash \bm{v} + (\bm{S} \circ \bm{O})\left[ - \bm{m} \oslash \left( (\bm{S} \circ \bm{O})^T \bm{v} \right) \circ (\mathbf{1}-\bm{B})  - \right. \\
\left. (r\alpha+1)\bm{m} \oslash \left( \bm{v}_0 + (\bm{S} \circ \bm{O})^T \bm{v} \right) \circ \bm{B} - \right. \\
\left. (r\alpha+1)\bm{L} \circ \left(\bm{I}^T \bm{v} \right) \oslash \left( \bm{v}_0 + (\bm{S} \circ \bm{O})^T \bm{v} \right)^{\circ 2} \right] + \\
\bm{I} \left[ - r(1-\alpha) \bm{m} \oslash \left( \bm{v}_0 + (\bm{S} \circ \bm{O})^T \bm{v} \right) \circ \bm{B} - \right. \\
\left. (r\alpha+1) \bm{L} \circ \left( (\bm{S} \circ \bm{O})^T \bm{v} \right) \oslash \left( \bm{v}_0 + (\bm{S} \circ \bm{O})^T \bm{v} \right)^{\circ 2} \right]
\ea
where $\oslash$, $\circ$, and $\circ 2$ denote the elementwise subtraction, addition, and square operations.

\newpage
\section*{Acknowledgments}
\noindent The authors would like to gratefully acknowledge Guillermo Gallego who suggested the idea of splitting the sales in case of partial availability.


\end{document}